\documentclass{amsart}
\usepackage{amssymb}

\input{prepictex}
\input{pictex}
\input{postpictex}

\theoremstyle{plain}
\newtheorem{theorem}{Theorem}[section]

\newtheorem{lemma}[theorem]{Lemma}
\newtheorem{prop}[theorem]{Proposition}
\theoremstyle{remark}
\newtheorem{notation}[theorem]{Notation}
\newtheorem{rmk}[theorem]{Remark}
\theoremstyle{definition}
\newtheorem{dfn}[theorem]{Definition}

\newtheorem{examples}[theorem]{Examples}

\numberwithin{equation}{section}

\DeclareMathOperator{\Obj}{Obj}
\DeclareMathOperator{\Mor}{Mor}
\DeclareMathOperator{\dom}{dom}
\DeclareMathOperator{\ran}{ran}
\DeclareMathOperator{\Ext}{Ext}
\DeclareMathOperator{\modulo}{mod}
\DeclareMathOperator{\Seg}{Seg}

\newcommand{\Lmin}{\Lambda^{\min}}

\newcommand{\field}[1]{\mathbb{#1}}
\newcommand{\NN}{\field{N}}

\newcommand{\ZZ}{\field{Z}}
\newcommand{\Bb}{{\mathcal B}}
\newcommand{\Gg}{{\mathcal G}}
\newcommand{\Oo}{{\mathcal O}}
\newcommand{\Uu}{{\mathcal U}}
\newcommand{\Zz}{{\mathcal Z}}
\newcommand{\CE}{{\mathcal{CE}}}
\newcommand{\Tt}{{\mathcal T}}

\begin{document}

\title[Groupoid models for the $C^*$-algebras of topological $k$-graphs]{Groupoid models for the $C^*$-algebras of topological higher-rank graphs}
\author{Trent Yeend}
\address{School of Mathematical \& Physical Sciences \\ University of Newcastle \\ NSW 2308 \\ AUSTRALIA}
\email{Trent.Yeend@newcastle.edu.au}


\begin{abstract}
We provide groupoid models for Toeplitz and Cuntz-Krieger algebras of topological higher-rank graphs. Extending the groupoid models used in the theory of graph algebras and topological dynamical systems to our setting, we prove results on essential freeness and amenability of the groupoids which capture the existing theory, and extend results involving group crossed products of graph algebras.
\end{abstract}

\maketitle

\section{Introduction}\label{sec:intro}

In recent years there has been significant interest in different generalizations of Cuntz-Krieger algebras, including the $C^*$-algebras of higher-rank graphs \cite{KP2,RSY1,RSY2,FMY1} and the $C^*$-algebras of topological graphs \cite{Kat1,Kat2}. In this article, we use groupoid models to explore a common approach to both generalizations.

We begin by introducing the notion of a topological higher-rank graph. Given a topological graph $\Lambda$ of rank $k\in\NN$, we define a topological path space $X_\Lambda$, which contains the finite paths of $\Lambda$ together with paths which are infinite in some or all of their $k$ dimensions. There is a natural action of $\Lambda$ on $X_\Lambda$ given by concatenation and removal of initial path segments, and using this we define a groupoid $G_\Lambda$, called the path groupoid of $\Lambda$, which has $X_\Lambda$ as its unit space.

We identify a topological analogue of the finitely aligned condition of \cite{RSY2,FMY1,PatWel}, called compactly aligned. If $\Lambda$ fails to be compactly aligned, then $G_\Lambda$ fails to be a topological groupoid: the range and source maps are not continuous and the topology on $G_\Lambda$ is not locally compact. However, if $\Lambda$ is compactly aligned, then $G_\Lambda$ is a locally compact topological groupoid which is $r$-discrete in the sense that the unit space $G_\Lambda^{(0)}$ is open in $G_\Lambda$. Furthermore, $G_\Lambda$ admits a Haar system, and so we may define the full groupoid $C^*$-algebra $C^*(G_\Lambda)$, which we refer to as the Toeplitz algebra of $\Lambda$.

We identify a closed invariant subset $\partial\Lambda$ of $G_\Lambda^{(0)}=X_\Lambda$, called the boundary-path space. The boundary-path groupoid of $\Lambda$ is the reduction $\Gg_\Lambda:= G_\Lambda|_{\partial\Lambda}$: a locally compact $r$-discrete groupoid admitting a Haar system. The Cuntz-Krieger algebra of $\Lambda$ is then defined to be the full groupoid $C^*$-algebra $C^*(\Gg_\Lambda)$. When $\Lambda$ is a finitely aligned discrete $k$-graph or the finite-path space of a second-countable topological graph, we recover the usual Toeplitz and Cuntz-Krieger algebras of the graph.

We then consider an analogue of the Aperiodicity Condition used in \cite{KP2}. We extend \cite[Proposition~4.5]{KP2} and \cite[Proposition~7.2]{FMY1} to our setting, proving that a compactly aligned topological $k$-graph $\Lambda$ satisfies the Aperiodicity Condition if and only if $\Gg_\Lambda$ is essentially free in the sense that the units with trivial isotropy are dense in $\Gg_\Lambda^{(0)}$.

Next, we address amenability of the boundary-path groupoid, showing that $\Gg_\Lambda$ is amenable if $\Lambda$ is either a finitely aligned discrete $k$-graph, a topological $1$-graph, or a \emph{proper} topological $k$-graph \emph{without sources}.

We end the article with a section on crossed products of topological $k$-graph algebras by coactions, extending \cite[Theorem~2.4]{KalQR} and \cite[Corollary~5.3]{KP2}. To begin, given a topological $k$-graph $\Lambda$, a locally compact group $A$ and a continuous functor $c:\Lambda\to A$, we define the notion of a skew-product topological $k$-graph $\Lambda\times_c A$.

If $A$ is abelian, then there are induced actions $\alpha$ of the dual group $\widehat{A}$ on $C^*(G_\Lambda)$ and $C^*(\Gg_\Lambda)$, and we extend \cite[Corollary~5.3]{KP2}, proving that the crossed product $C^*$-algebras $C^*(G_\Lambda)\times_\alpha \widehat{A}$ and $C^*(\Gg_\Lambda)\times_\alpha \widehat{A}$ are isomorphic to $C^*(G_{\Lambda\times_c A})$ and $C^*(\Gg_{\Lambda\times_c A})$, respectively.

If $A$ is discrete, then there are induced coactions $\delta$ of $A$ on $C^*(G_\Lambda)$ and $C^*(\Gg_\Lambda)$, and we extend \cite[Theorem~2.4]{KalQR}, proving that $C^*(G_\Lambda)\times_\delta A\cong C^*(G_{\Lambda\times_c A})$ and $C^*(\Gg_\Lambda)\times_\delta A\cong C^*(\Gg_{\Lambda\times_c A})$.

\section{Topological higher-rank graphs}\label{sec:top h-r graphs}

\begin{dfn}
Given $k\in\NN$, a \emph{topological $k$-graph} is a pair $(\Lambda,d)$ consisting of a small category $\Lambda=(\Obj(\Lambda),\Mor(\Lambda),r,s)$ and a functor $d:\Lambda\to\NN^k$, called the \emph{degree map}, which satisfy the following:
\begin{enumerate}
\item
$\Obj(\Lambda)$ and $\Mor(\Lambda)$ are second-countable locally compact Hausdorff spaces;
\item
$r,s:\Mor(\Lambda)\to\Obj(\Lambda)$ are continuous and $s$ is a local homeomorphism;
\item
Composition $\circ : \Lambda\times_c\Lambda\to\Lambda$ is continuous and open, where $\Lambda\times_c\Lambda$ has the relative topology inherited from the product topology on $\Lambda\times\Lambda$;
\item
$d$ is continuous, where $\NN^k$ has the discrete topology;
\item
For all $\lambda\in\Lambda$ and $m,n\in\NN^k$ such that $d(\lambda)=m+n$, there exists unique $(\xi,\eta)\in\Lambda\times_c\Lambda$ such that $\lambda=\xi\eta$, $d(\xi)=m$ and $d(\eta)=n$.
\end{enumerate}
We refer to the morphisms of $\Lambda$ as \emph{paths} and to the objects of $\Lambda$ as \emph{vertices}. The codomain and domain maps from $\Lambda$ are called the range and source maps, respectively.
\end{dfn}

\begin{notation}
For $m\in\NN^k$, we write $m_i$ for the $i$th coordinate of $m$. We use the partial ordering $\le$ on $\NN^k$ defined by $m\le n \iff m_i\le n_i$ for all $i\in\{1,\dots,k\}$, so least upper bounds and greatest lower bounds are given by $(m\vee n)_i=\max\{m_i,n_i\}$ and $(m\wedge n)_i = \min\{m_i,n_i\}$, respectively.
For $m\in\NN^k$, define $\Lambda^m$ to be the set $d^{-1}(\{m\})$ of paths of degree $m$. Define $\Lambda *_s\Lambda :=\{(\lambda,\mu)\in\Lambda\times\Lambda : s(\lambda)=s(\mu)\}$, and for $U,V\subset\Lambda$ define $U*_sV:=(U\times V)\cap (\Lambda *_s\Lambda)$. For $p,q\in\NN^k$, $U\subset\Lambda^p$ and $V\subset\Lambda^q$, we write
\[
U\vee V := U\Lambda^{(p\vee q)-p}\cap V\Lambda^{(p\vee q)-q}
\]
for the set of \emph{minimal common extensions} of paths from $U$ and $V$. For $\lambda,\mu\in\Lambda$, we write
\[
\Lmin(\lambda,\mu):=\{(\alpha,\beta) | \lambda\alpha=\mu\beta, d(\lambda\alpha)=d(\lambda)\vee d(\mu)\}
\]
for the set of pairs which give minimal common extensions of $\lambda$ and $\mu$; that is,
\[
\Lmin(\lambda,\mu) = \{(\alpha,\beta) | \lambda\alpha=\mu\beta\in \{\lambda\}\vee\{\mu\}\}.
\]
\end{notation}

\begin{dfn}
A topological $k$-graph $(\Lambda,d)$ is \emph{compactly aligned} if for all $p,q\in\NN^k$ and compact $U\subset\Lambda^p$ and $V\subset\Lambda^q$, the set $U\vee V$ is compact.
\end{dfn}

\begin{rmk}
A discrete $k$-graph $(\Lambda,d)$ is \emph{finitely aligned} if for all $\lambda,\mu\in\Lambda$, the set $\Lmin(\lambda,\mu)$ is finite; these discrete $k$-graphs form the scope of the higher-rank graph theory to date (see \cite{RSY2,PatWel,FMY1}). As compactness is equivalent to finiteness for discrete topologies, it follows that a discrete $k$-graph is compactly aligned if and only if it is finitely aligned.
\end{rmk}

\begin{examples}\label{ex:graph examples}
\begin{enumerate}
\item\label{ex:item:discrete k-graph}
Any higher-rank graph (as defined in \cite{KP2,RSY2} etc.) may be regarded as a topological higher-rank graph with discrete topologies on the object and morphism sets.
\item\label{ex:item:top graph}
Let $E$ be a second-countable topological graph as defined in \cite[Definition~2.1]{Kat1}; that is, $E=(E^0,E^1,r,s)$ is a directed graph with $E^0, E^1$ second-countable locally compact Hausdorff spaces, $r,s:E^1\to E^0$ continuous, and $s$ a local homeomorphism. The free category generated by $E$, endowed with the relative topology inherited from the union of the product topologies, together with the length function $l(e_1\cdots e_n) = n$, forms a topological $1$-graph $(E^*,l)$. Conversely, given a topological $1$-graph $(\Lambda,d)$, the quadruple $E_\Lambda := (\Lambda^0,\Lambda^1,r|_{\Lambda^1},s|_{\Lambda^1})$ is a second-countable topological graph with $((E_\Lambda)^*,l) \cong (\Lambda,d)$.
\item\label{ex:item:SGDS}
Let $(X,\theta)$ be a singly generated dynamical system as defined in \cite{ArRen,Ren4}; that is, $X$ is a second-countable locally compact Hausdorff space and $\theta$ is a local homeomorphism from an open subset $\dom(\theta)$ of $X$ onto an open subset $\ran(\theta)$ of $X$. In \cite[Section~10.3]{Kat1}, Katsura constructs a topological graph $E(X,\theta)$ by setting $E(X,\theta)^0 := X$ and $E(X,\theta)^1 := \dom(\theta)$, and for $x\in E(X,\theta)^1$, setting $r(x) := x$ and $s(x) := \theta(x)$. So, as was done for a general topological graph in the previous example, we may form the topological $1$-graph $(\Lambda(X,\theta),d) := (E(X,\theta)^*,l)$. The following example generalizes this construction.
\item\label{ex:item:k-GDS}
Let $X$ be a second-countable locally compact Hausdorff space, and for $i=1,\dots,k$, let $\theta_i$ be a local homeomorphism from an open subset $\dom(\theta_i)$ of $X$ onto an open subset $\ran(\theta_i)$ of $X$, such that for all $i,j\in \{1,\dots,k\}$,
\[
\dom(\theta_j\theta_i)=\theta_i^{-1}(\ran(\theta_i)\cap\dom(\theta_j)) = \theta_j^{-1}(\ran(\theta_j)\cap \dom(\theta_i))=\dom(\theta_i\theta_j)
\]
and for $x\in \dom(\theta_j\theta_i)$,
\[
\theta_j(\theta_i(x))=\theta_i(\theta_j(x)).
\]
We define a topological $k$-graph $(\Lambda(X,\{\theta_i\}_{i=1}^k),d)$ by setting
\[
\Obj(\Lambda(X,\{\theta_i\})) := X,
\]
\begin{align*}
\Mor(\Lambda(X,\{\theta_i\})) &:= \{(m,x)\in\NN^k\times X | x\in\dom(\theta_1^{m_1}\theta_2^{m_2}\cdots\theta_k^{m_k})\} \\
&= \bigcup_{m\in\NN^k} \{m\}\times\dom(\theta_1^{m_1}\cdots\theta_k^{m_k}),
\end{align*}
\[
r(m,x):= x,\quad s(m,x) := \theta_1^{m_1}\cdots\theta_k^{m_k}(x)
\]
and
\[
(m,x)(n,\theta_1^{m_1}\cdots\theta_k^{m_k}(x)) := (m+n,x),
\]
giving $\Lambda(X,\{\theta_i\}_{i=1}^k)$ the relative topology inherited from the product topology, and setting $d(m,x):= m$.
\end{enumerate}
\end{examples}

\section{The path groupoid}\label{sec:path groupoid}

We begin this section by associating a groupoid $G_\Lambda$ to each topological $k$-graph $(\Lambda,d)$. We first define the unit space $G_\Lambda^{(0)}$ as a space $X_\Lambda$ of paths of $\Lambda$; the finite paths in $X_\Lambda$ are characterized by the morphisms $\lambda\in\Lambda$, however we must also consider paths of $\Lambda$ which are infinite in some or all of their $k$-dimensions. To do this, we first define appropriate rank-$k$ path prototypes for each and every degree -- finite, infinite and partially infinite -- and then obtain $X_\Lambda$ as the set of representations of the path prototypes. The morphisms of the category $\Lambda$ are then in correspondence with the representations of those path prototypes which are finite in each of the $k$ dimensions.

For $k\in\NN$ and $m\in (\NN\cup\{\infty\})^k$, define the topological $k$-graph $(\Omega_{k,m},d)$ by giving the discrete topologies to the sets
\[
\Obj(\Omega_{k,m}):=\{p\in\NN^k | p\le m\}
\]
and
\[
\Mor(\Omega_{k,m}):= \{(p,q)\in\NN^k\times\NN^k | p\le q\le m\},
\]
and setting $r(p,q):=p$, $s(p,q):=q$, $(n,p)\circ (p,q):= (n,q)$ and $d(p,q):= q-p$.

Let $(\Lambda_1,d_1)$ and $(\Lambda_2,d_2)$ be topological $k$-graphs.  A \emph{graph morphism} between $\Lambda_1$ and $\Lambda_2$ is a continuous functor $x:\Lambda_1\to\Lambda_2$ satisfying $d_2(f(\lambda))=d_1(\lambda)$ for all $\lambda\in\Lambda_1$.

\begin{dfn}\label{dfn:defining X}
Let $(\Lambda,d)$ be a topological $k$-graph. We define
\[
X_\Lambda :=\bigcup_{m\in(\NN\cup\{\infty\})^k}\{x:\Omega_{k,m}\to\Lambda ~|~ x \text{ is a graph morphism}\}.
\]
We extend the range and degree maps to $x:\Omega_{k,m}\to\Lambda$ in $X_\Lambda$ by setting $r(x):=x(0)$ and $d(x):=m$. For $v\in\Lambda^0$ we define $v X_\Lambda:=\{x\in X_\Lambda ~|~ r(x)=v\}$.
\end{dfn}

\begin{notation}
For each $\lambda\in\Lambda$ there is a unique graph morphism $x_\lambda:\Omega_{k,d(\lambda)}\to\Lambda$ such that $x_\lambda(0,d(\lambda))=\lambda$; in this sense, we may view $\Lambda$ as a subset of $X_\Lambda$, and we refer to elements of $X_\Lambda$ as paths. Indeed, for $\lambda\in\Lambda$ and $p,q\in\NN^k$ with $0\le p\le q\le d(\lambda)$, we may write $\lambda(0,p)$, $\lambda(p,q)$ and $\lambda(q,d(\lambda))$ for the unique elements of $\Lambda$ which satisfy $\lambda=\lambda(0,p)\lambda(p,q)\lambda(q,d(\lambda))$, $d(\lambda(0,p))=p$, $d(\lambda(p,q))=q-p$ and $d(\lambda(q,d(\lambda)))=d(\lambda)-q$.
\end{notation}

Straightforward arguments give the following lemma.

\begin{lemma}
Let $(\Lambda,d)$ be a topological $k$-graph. For $x\in X_\Lambda$, $m\in\NN^k$ with $m\le d(x)$, and $\lambda\in\Lambda$ with $s(\lambda)=r(x)$, there exist unique paths $\lambda x$ and $\sigma^m x$ in $X_\Lambda$ satisfying $d(\lambda x)=d(\lambda)+d(x)$, $d(\sigma^m x) = d(x)-m$,
\[
(\lambda x)(0,p) = \begin{cases}
\lambda(0,p) &\text{if } p\le d(\lambda) \\
\lambda x(0,p-d(\lambda)) &\text{if } d(\lambda)\le p\le d(\lambda x),
\end{cases}
\]
and
\[
(\sigma^m x)(0,p) = x(m,m+p) \quad\text{for } p\le d(\sigma^m x).
\]
\end{lemma}

\begin{dfn}
Let $(\Lambda,d)$ be a topological $k$-graph. Define the \emph{path groupoid} $G_\Lambda$ to be the groupoid with object set $\Obj(G_\Lambda):= X_\Lambda$, morphism set
\begin{align*}
\Mor(G_\Lambda) &:= \{ (\lambda x, d(\lambda)-d(\mu), \mu x) \in X_\Lambda\times\ZZ^k\times X_\Lambda ~|~ \\
&\phantom{MMMMMMNMMN}(\lambda,\mu)\in\Lambda *_s\Lambda,~ x\in X_\Lambda \text{ and } s(\lambda)=r(x)\} \\
&=\{ (x,m,y) \in X_\Lambda\times\ZZ^k\times X_\Lambda ~|~ \text{ there exist } p,q\in\NN^k \text{ such that}\\
&\phantom{MMMNNMMN} p\le d(x),~ q\le d(y),~ p-q=m \text{ and } \sigma^p x = \sigma^q y \},
\end{align*}
range and source maps
\[
r(x,m,y):= x \text{ and } s(x,m,y):= y,
\]
composition
\[
((x,m,y),(y,n,z))\mapsto (x,m+n,z),
\]
and inversion
\[
(x,m,y)\mapsto (y,-m,x).
\]
\end{dfn}

\begin{notation}
Let $(\Lambda,d)$ be a topological $k$-graph. For $F\subset \Lambda *_s\Lambda$ and $m\in\ZZ^k$, define $Z(F,m)\subset G_\Lambda$ by
\[
Z(F,m) := \{(\lambda x,d(\lambda)-d(\mu),\mu x)\in G_\Lambda ~|~ (\lambda,\mu)\in F, d(\lambda)-d(\mu)=m\}.
\]
For $U\subset\Lambda$, define $Z(U)\subset G_\Lambda^{(0)}$ by
\[
Z(U) := Z(U*_sU,0)\cap Z(\Lambda^0*_s\Lambda^0,0).
\]
\end{notation}

\begin{prop}\label{prop:defining a topology on G_Lambda}
Let $(\Lambda,d)$ be a topological $k$-graph. The family of sets of the form
\[
Z(U*_s V,m)\cap Z(F,m)^c,
\]
where $m\in\ZZ^k$, $U,V\subset\Lambda$ are open and $F\subset \Lambda *_s\Lambda$ is compact, is a basis for a second-countable Hausdorff topology on $G_\Lambda$.
\end{prop}

\begin{proof}
To see that the family of sets forms a basis, suppose
\[
(x,m,y)\in \big( Z(U_1*_sV_1,m)\cap Z(F_1,m)^c\big)\cap\big( Z(U_2*_sV_2,m)\cap Z(F_2,m)^c\big),
\]
where $m\in\ZZ^k$, $U_1,U_2,V_1,V_2\subset\Lambda$ are open and $F_1,F_2\subset\Lambda*_s\Lambda$ are compact. We then have the existence of $(\lambda,\mu)\in U_1*_sV_1$, $(\xi,\eta)\in U_2*_sV_2$ and $w,z\in X_\Lambda$ such that
\[
(x,m,y) = (\lambda w,d(\lambda)-d(\mu),\mu w) = (\xi z,d(\xi)-d(\eta),\eta z).
\]
Hence the pair $w,z$ extend $\lambda,\xi$ and $\mu,\eta$ to common paths $x$ and $y$, respectively, and setting
\[
\alpha:= w(0,(d(\lambda)\vee d(\xi))-d(\lambda)) = w(0,(d(\mu)\vee d(\eta))-d(\mu))
\]
and
\[
\beta:= z(0,(d(\lambda)\vee d(\xi))-d(\xi)) = z(0,(d(\mu)\vee d(\eta))-d(\eta)),
\]
we have
\[
\lambda\alpha=\xi\beta \quad\text{and}\quad \mu\alpha=\eta\beta.
\]
Let $W_1\subset\Lambda$ be an open neighbourhood of $\alpha$ such that $s|_{W_1}$ is a homeomorphism, and let $W_2\subset\Lambda$ be an open neighbourhood of $\beta$ such that $s|_{W_2}$ is a homeomorphism. Since composition is open, the sets $U_1W_1, V_1W_1, U_2W_2, V_2W_2\subset\Lambda$ are open, and we have
\[
(x,m,y)\in Z\big( (U_1W_1*_sV_1W_1)\cap (U_2W_2*_sV_2W_2),m\big)\cap Z(F_1\cup F_2,m)^c.
\]
Furthermore, since $s|_{W_1}$ and $s|_{W_2}$ are homeomorphisms, it follows that
\begin{align*}
Z\big( (U_1W_1*_s&V_1W_1)\cap (U_2W_2*_sV_2W_2),m\big)\cap Z(F_1\cup F_2,m)^c \\
&\subset \big( Z(U_1*_sV_1,m)\cap Z(F_1,m)^c\big)\cap\big( Z(U_2*_sV_2)\cap Z(F_2,m)^c\big),
\end{align*}
as required.

Second-countability is clear. It remains to show the topology is Hausdorff. Let $(w,m,x)$ and $(y,n,z)$ be distinct elements of $G_\Lambda$. If $m\neq n$, then $(w,m,x)\in Z(\Lambda *_s\Lambda,m)$, $(y,n,z)\in Z(\Lambda *_s\Lambda,n)$ and $Z(\Lambda *_s\Lambda,m)\cap Z(\Lambda *_s\Lambda,m) = \emptyset$, so we assume $m=n$. Furthermore, if $r(w)\neq r(y)$, then taking open neighbourhoods $U,V\subset\Lambda^0$ of $r(w)$ and $r(y)$, respectively, such that $U\cap V = \emptyset$, we have $(w,m,x)\in Z(U\Lambda*_s\Lambda,m)$, $(y,m,z)\in Z(V\Lambda*_s\Lambda,m)$ and $Z(U\Lambda*_s\Lambda,m)\cap Z(V\Lambda *_s\Lambda,m)=\emptyset$. A similar argument holds if $r(x)\neq r(z)$, so we assume $r(w)=r(y)$ and $r(x)=r(z)$.

We must have either $w\neq y$ or $x\neq z$, so assume $w\neq y$. Let $p\in\NN^k$ be minimal with respect to the conditions
\begin{equation*}
\begin{split}
&p\le d(w)\wedge d(y),~ w(0,p)=y(0,p), \text{ and } \\
&w(0,p+e_i)\neq y(0,p+e_i) \text{ for some } i\in\{1,\dots,k\}.
\end{split}
\end{equation*}
We must have either $d(w)\ge p+e_i$ or $d(y)\ge p+e_i$: If both $d(w)\ge p+e_i$ and $d(y)\ge p+e_i$, we can take open neighbourhoods $U,V\subset\Lambda^{p+e_i}$ of $w(0,p+e_i)$ and $y(0,p+e_i)$, respectively, such that $U\cap V = \emptyset$. Then $(w,m,x)\in Z(U\Lambda*_s\Lambda,m)$, $(y,m,z)\in Z(V\Lambda*_s\Lambda,m)$ and $Z(U\Lambda*_s\Lambda,m)\cap Z(V\Lambda*_s\Lambda,m)=\emptyset$, so we assume $d(w)\ge p+e_i$ and $d(y)\not\ge p+e_i$.

Express $(w,m,x)$ as $(\lambda w', d(\lambda)-d(\mu),\mu w')$, where $d(\lambda)\ge p+e_i$, and let $U\subset\Lambda^{d(\lambda)}$ and $V\subset\Lambda^{d(\mu)}$ be relatively compact open neighbourhoods of $\lambda$ and $\mu$, respectively. Then $(w,m,x)\in Z(U*_sV,m)$, $(y,m,z)\in Z(\overline{U}*_s\overline{V},m)^c$ and $Z(U*_sV,m)\cap Z(\overline{U}*_s\overline{V},m)^c=\emptyset$, proving the topology is Hausdorff.
\end{proof}

\begin{rmk}\label{rmk:not compactly aligned entails problems}
If $(\Lambda,d)$ is not compactly aligned, then the topology on $G_\Lambda$ defined by Proposition~\ref{prop:defining a topology on G_Lambda} may not be locally compact, and, under this topology, $G_\Lambda$ may not be a topological groupoid. To illustrate these two facts, we consider two $2$-graphs which fail to be compactly aligned. We describe the $2$-graphs in terms of their $1$-skeletons as in \cite[Section~2]{RSY1}.

Let $(\Lambda_1,d)$ be the discrete topological $2$-graph with $1$-skeleton:
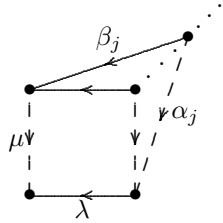
\[
\beginpicture
\put {$\bullet$} at 0 0
\put {$\bullet$} at 40 0
\put {$\bullet$} at 0 40
\put {$\bullet$} at 40 40
\put {$\bullet$} at 60 60
\put {$\lambda$} at 20 -5
\put {$\mu$} at -5 20
\put {$\alpha_j$} at 59 30
\put {$\beta_j$} at 30 59
\put {$.$} at 45 45
\put {$.$} at 50 50
\put {$.$} at 55 55
\put {$.$} at 64 64
\put {$.$} at 68 68
\put {$.$} at 72 72
\arrow <0.15cm> [0.25,0.75] from  0 20 to 0 19.5
\arrow <0.15cm> [0.25,0.75] from  40 20 to 40 19.5
\arrow <0.15cm> [0.25,0.75] from  20 0 to 19.5 0
\arrow <0.15cm> [0.25,0.75] from  20 40 to 19.5 40
\arrow <0.15cm> [0.25,0.75] from  50 30 to 49.5 28.5
\arrow <0.15cm> [0.25,0.75] from  30 50 to 28.5 49.5
\plot 0 0 40 0 /
\plot 0 40 40 40 /
\plot 0 40 60 60 /
\setdashes
\plot 0 0 0 40 /
\plot 40 0 40 40 /
\plot 40 0 60 60 /
\endpicture
\]
So $\lambda\alpha_j=\mu\beta_j$ for all $j\in\NN$. The sequence $\langle (\lambda\alpha_j,(1,0),\alpha_j)\rangle_{j\in\NN}$ converges to $(\lambda,(1,0),s(\lambda))$ in $G_{\Lambda_1}$, but $\langle r(\lambda\alpha_j,(1,0),\alpha_j)\rangle_{j\in\NN} = \langle \lambda\alpha_j\rangle_{j\in\NN}$ does not converge to $r(\lambda,(1,0),s(\lambda))=\lambda$ in $G_{\Lambda_1}^{(0)}$ since $\{\lambda\}= Z(\lambda)\cap Z(\mu)^c$. Therefore the range map in $G_{\Lambda_1}$ is not continuous.

Now, taking a family of copies of $\Lambda_1$ indexed by $\NN$, and identifying the path $\lambda$ from each, we obtain the $2$-graph $(\Lambda_2,d)$ with $1$-skeleton:
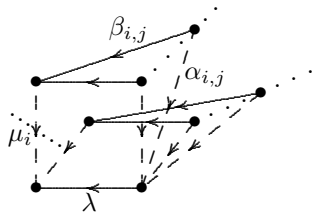
\[
\beginpicture
\put {$\bullet$} at 0 0
\put {$\bullet$} at 40 0
\put {$\bullet$} at 20 25
\put {$\bullet$} at 60 25
\put {$\bullet$} at 85 36
\put {$\bullet$} at 0 40
\put {$\bullet$} at 40 40
\put {$\bullet$} at 60 60
\put {$\lambda$} at 20 -5
\put {$\mu_i$} at -6 19
\put {$\beta_{i,j}$} at 35 60
\put {$\alpha_{i,j}$} at 64 42
\put {$.$} at 45 45
\put {$.$} at 50 50
\put {$.$} at 55 55
\put {$.$} at 64 64
\put {$.$} at 68 68
\put {$.$} at 72 72
\put {$.$} at 3 20
\put {$.$} at 6 18
\put {$.$} at 9 16
\put {$.$} at -3 24
\put {$.$} at -6 26
\put {$.$} at -9 28
\put {$.$} at 66.25 27.75
\put {$.$} at 72.5 30.5
\put {$.$} at 78.75 33.25
\put {$.$} at 91.25 38.75
\put {$.$} at 97.5 41.5
\put {$.$} at 103.75 44.25
\arrow <0.15cm> [0.25,0.75] from  0 20 to 0 19.5
\arrow <0.15cm> [0.25,0.75] from  40 20 to 40 19.5
\arrow <0.15cm> [0.25,0.75] from  20 0 to 19.5 0
\arrow <0.15cm> [0.25,0.75] from  20 40 to 19.5 40
\arrow <0.15cm> [0.25,0.75] from  50.5 31.5 to 50 30
\arrow <0.15cm> [0.25,0.75] from  30 50 to 28.5 49.5
\arrow <0.15cm> [0.25,0.75] from  12 15 to 11 13.75
\arrow <0.15cm> [0.25,0.75] from  53.5 30.670 to 53 30.585
\arrow <0.15cm> [0.25,0.75] from  42.5 25 to 42 25
\arrow <0.15cm> [0.25,0.75] from  50.5 13.125 to 50 12.5
\arrow <0.15cm> [0.25,0.75] from  58.5 14.8 to 58 14.4
\plot 60 25 20 25 /
\plot 85 36 20 25 /
\plot 0 0 40 0 /
\plot 0 40 40 40 /
\plot 0 40 60 60 /
\setdashes
\plot 20 25 0 0 /
\plot 60 25 40 0 /
\plot 85 36 40 0 /
\plot 0 0 0 40 /
\plot 40 0 40 40 /
\plot 40 0 60 60 /
\endpicture
\]
So, $\lambda\alpha_{i,j}=\mu_i\beta_{i,j}$ for all $i,j\in\NN$. We claim that the unit $\lambda\in G_{\Lambda_2}^{(0)}$ has no compact neighbourhood: First note that any neighbourhood of $\lambda$ in $G_{\Lambda_2}^{(0)}$ contains a basis set of the form
\[
Z(\lambda)\cap\bigcap_{k=1}^nZ(\mu_{i_k})^c,\quad\text{where } i_k\in\NN \text{ for } k=1,\dots,n;
\]
furthermore, given such a basis set, and choosing $l\in\NN$ with $l\neq i_k$ for all $k\in\{1,\dots,n\}$, the family
\[
\Big\{ Z(\lambda)\cap\bigcap_{k=1}^nZ(\mu_{i_k})^c\cap Z(\mu_l)^c\Big\}\cup \big\{ Z(\mu_l\beta_{l,j}) ~|~ j\in\NN\big\}
\]
forms an infinite, disjoint, open cover which, consequently, has no finite subcover. Hence the topology on $G_{\Lambda_2}$ is not locally compact.
\end{rmk}

The following two lemmas allow us to restrict the type of basis elements we will need to consider. We omit the proof of the first lemma.

\begin{lemma}[{cf. \cite[Remark~5.5]{FMY1}}] \label{lem:choosing better basis sets for unit space}
Let $(\Lambda,d)$ be a compactly aligned topological $k$-graph. For any relatively compact $U\subset\Lambda$ and compact $F\subset\Lambda$,
\[
Z(U)\cap Z(F)^c = Z(U)\cap Z(\overline{U}\vee F)^c.
\]
Thus, since $\Lambda$ is locally compact and since the source map in $\Lambda$ is a local homeomorphism, we need only consider basis sets for $G_\Lambda^{(0)}$ of the form $Z(U)\cap Z(F)^c$, where $U\subset\Lambda$ is relatively compact and open, and $F\subset\Lambda$ is compact and satisfies $\mu\in F$ implies $\mu = \lambda\mu'$ for some $\lambda\in\overline{U}$.
\end{lemma}

\begin{lemma}\label{lem:choosing better basis sets}
Let $(\Lambda,d)$ be a compactly aligned topological $k$-graph. Let $p,q\in\NN^k$, let $U\subset\Lambda^p$, $V\subset\Lambda^q$ be relatively compact open sets, and let $F\subset\Lambda *_s\Lambda$ be compact. There exists a compact set $F'\subset\Lambda *_s\Lambda$ such that
\[
Z(U*_sV,p-q)\cap Z(F,p-q)^c = Z(U*_sV,p-q)\cap Z(F',p-q)^c
\]
and
\[
(\xi,\eta)\in F' \text{ implies } (\xi,\eta)=(\lambda\alpha,\mu\alpha) \text{ for some } (\lambda,\mu)\in \overline{U}*_s\overline{V} \text{ and } \alpha\in\Lambda.
\]
\end{lemma}

\begin{proof}
We can assume $d(\xi)-d(\eta)=m$ for all $(\xi,\eta)\in F$. Since $F$ is compact, there exist $m^1,\dots,m^l\in\NN^k$ such that each $m^j\ge p-q$ and
\[
F = \bigcup_{j=1}^l F\cap (\Lambda^{m^j}*_s\Lambda^{m^j-(p-q)});
\]
for $j=1,\dots,l$, define $F_j:= F\cap (\Lambda^{m^j}*_s\Lambda^{m^j-(p-q)})$ and
\[
\begin{split}
F_j':= \{(\xi,\eta)&\in \Lambda^{m^j\vee p}*_s\Lambda^{(m^j-p+q)\vee q} ~|~ (\xi(0,m^j),\eta(0,m^j-p+q))\in F_j, \\
&(\xi(0,p),\eta(0,q))\in \overline{U}*_s\overline{V} \text{ and } \xi(p,m^j\vee p)=\eta(q,(m^j-p+q)\vee q)\}.
\end{split}
\]
Letting $P_1,P_2:\Lambda*_s\Lambda\to\Lambda$ be the coordinate projections, we see that each $F_j'$ is a closed subset of the compact set $(P_1(F_j)\vee\overline{U})*_s(P_2(F_j)\vee\overline{V})$. Hence each $F_j'$ is compact. Defining $F':=\bigcup_{j=1}^lF_j'$ completes the proof.
\end{proof}

\begin{dfn}
Let $(\Lambda,d)$ be a topological $k$-graph. For $m,p,q\in\NN^k$ with $p\le q\le m$, define the continuous map $\Seg^m_{(p,q)}:\Lambda^m\to\Lambda^{q-p}$ by $\Seg^m_{(p,q)}(\lambda):=\lambda(p,q)$.
\end{dfn}

\begin{dfn}
Let $(\Lambda,d)$ be a topological $k$-graph. An infinite sequence of paths in $\Lambda$ is \emph{wandering} if for every compact set $F\subset\Lambda$, the sequence is eventually in $\Lambda\setminus F$; that is to say, the sequence visits any compact set at most finitely many times.
\end{dfn}

We have the following technical characterization of convergence in $G_\Lambda^{(0)}$ (cf. \cite[Remark~5.6]{FMY1} and \cite[page~653]{Pat2}).

\begin{prop}
Let $(\Lambda,d)$ be a compactly aligned topological $k$-graph, and let $\langle x_j\rangle_{j\in\NN}$ and $x$ be in $G_\Lambda^{(0)}$. Then
\[
\lim_{j\to\infty} x_j=x
\]
if and only if the following two conditions hold:
\begin{enumerate}
\item
for all $m\in\NN^k$ with $m\le d(x)$,
\[
\lim_{j\to\infty}x_j(0,m\wedge d(x_j)) = x(0,m);
\]
\item
for all $i\in\{1,\dots,k\}$ with $d(x)_i<\infty$ and for all $n\in\NN^k$ with $n\le d(x)$ and $n_i=d(x)_i$, if the set
\[
J(n,i):= \{j\in\NN ~|~ d(x_j)\ge n+e_i\}
\]
is infinite, then $\langle x_j(n,n+e_i)\rangle_{j\in J(n,i)}$ is wandering.
\end{enumerate}
\end{prop}

\begin{proof}
Assume $\displaystyle{\lim_{j\to\infty} x_j = x}$. For any $m\in\NN^k$ with $m\le d(x)$ and any open neighbourhood $U$ of $x(0,m)$, we have $x\in Z(U)$. Hence $\langle x_j\rangle_{j\in\NN}$ is eventually in $Z(U)$, it follows that $\langle x_j(0,m\wedge d(x_j)\rangle_{j\in\NN}$ is eventually in $U$, and Condition~(1) holds.

To show that Condition~(2) holds, suppose $d(x)<\infty$ for some $i\in\{1,\dots,k\}$, let $n\in\NN^k$ satisfy $n\le d(x)$ and $n_i=d(x)_i$, suppose $J(n,i) = \{j\in\NN ~|~ d(x_j)\ge n+e_i\}$ is infinite, and let $W\subset\Lambda^{e_i}$ be compact; we show that $\langle x_j(n,n+e_i)\rangle_{j\in J(n,i)}$ is eventually in $\Lambda^{e_i}\setminus W$.

Let $U\subset\Lambda^n$ be a relatively compact open neighbourhood of $x(0,n)$. Then $x\in Z(U)\cap Z(\overline{U}W)^c$, and hence, eventually, so is $\langle x_j\rangle_{j\in\NN}$. As $\langle x_j(0,n)\rangle_{j\in\NN}$ is eventually in $U$, it follows that $\langle x_j(n,n+e_i)\rangle_{j\in J(n,i)}$ is eventually in $\Lambda^{e_i}\setminus W$, as required.

Conversely, assume Conditions~(1) and~(2) hold. We argue by contradiction, so suppose there exist $m\in\NN^k$, relatively compact open $U\subset\Lambda^m$ and compact $F\subset\Lambda$ such that
\begin{equation}\label{eqn:x in nbd}
x\in Z(U)\cap Z(F)^c
\end{equation}
and
\begin{equation}\label{eqn:x_j sequence not in nbd}
x_j\not\in Z(U)\cap Z(F)^c \quad\text{for infinitely many } j\in\NN.
\end{equation}
By Condition~(1), we must have $x_j\in Z(U)$ eventually, so it follows that
\begin{equation}\label{eqn:x_j sequence in compact set}
x_j\in Z(F) \quad\text{for infinitely many } j\in\NN.
\end{equation}
Since $F$ is compact, without loss of generality we assume $F\subset \Lambda^M$ for some $M\in\NN^k$; by Lemma~\ref{lem:choosing better basis sets for unit space}, we can assume $m\le M$, retaining \eqref{eqn:x in nbd} through \eqref{eqn:x_j sequence in compact set}\footnote{This is the only point in the proof which relies on $(\Lambda,d)$ being compactly aligned.}.

We claim that $M\not\le d(x)$; otherwise, \eqref{eqn:x in nbd} implies $x(0,M)\in\Lambda^M\setminus F$, which combines with \eqref{eqn:x_j sequence in compact set} to contradict Condition~(1). Thus $M\not\le d(x)$.

Define $N:=M\wedge d(x)$. Then $N\le d(x)$ and there exists $i\in\{1,\dots,k\}$ such that $N_i=d(x)_i$ and $N+e_i\le M$. We have that $J(N,i) = \{j\in\NN ~|~ d(x_j)\ge N+e_i\}$ is infinite since it contains the infinite set $\{j\in\NN ~|~ x_j\in Z(F)\}$. By \eqref{eqn:x_j sequence in compact set}, we also have $x_j(N,N+e_i)$ in the compact set $\Seg^M_{(N,N+e_i)}(F)$ for infinitely many $j\in J(N,i)$, which contradicts Condition~(2). Therefore $\langle x_j\rangle_{j\in\NN}$ is eventually in $Z(U)\cap Z(F)^c$, and $\displaystyle{\lim_{j\to\infty}x_j=x}$.
\end{proof}

We deduce the following characterization of convergence in $G_\Lambda$.

\begin{prop}\label{prop:convergence means}
Let $(\Lambda,d)$ be a compactly aligned topological $k$-graph, let $p,q\in\NN^k$, and let the sequence $\langle (x_j,p-q,y_j)\rangle_{j\in\NN}$ and point $(x,p-q,y)$ be contained in $Z(\Lambda^p*_s\Lambda^q,p-q)$. Then
\[
\lim_{j\to\infty}(x_j,p-q,y_j) = (x,p-q,y)
\]
if and only if the following two conditions hold:
\begin{enumerate}
\item
for all $m\in\NN^k$,
\begin{enumerate}
\item
$\displaystyle{\lim_{j\to\infty} x_j(0,m\wedge d(x)\wedge d(x_j)) = x(0,m\wedge d(x))}$ and
\item
$\displaystyle{\lim_{j\to\infty} y_j(0,m\wedge d(y)\wedge d(y_j)) = y(0,m\wedge d(y))}$;
\end{enumerate}
\item
for all $i\in\{1,\dots, k\}$ with $d(x)_i<\infty$ and for all $n\in\NN^k$ with $p\le n\le d(x)$ and $n_i=d(x)_i$, if the set
\[
J(n,i):=\{ j\in\NN ~|~ d(x_j)\ge n+e_i\}
\]
is infinite, then
\[
\langle x_j(n,n+e_i)\rangle_{j\in J(n,i)}
\]
is wandering.
\end{enumerate}
\end{prop}

\begin{rmk}
Condition~(2) of Proposition~\ref{prop:convergence means} is equivalent to the following condition stated in terms of $y$ and the $y_j$:\emph{
\begin{itemize}
\item[$(2')$]
for all $i\in\{1,\dots, k\}$ with $d(y)_i<\infty$ and for all $n\in\NN^k$ with $q\le n\le d(y)$ and $n_i=d(y)_i$, if the set
\[
J(n,i)=\{j\in\NN ~|~ d(y_j)\ge n+e_i\}
\]
is infinite, then
\[
\langle y_j(n,n+e_i)\rangle_{j\in J(n,i)}
\]
is wandering.
\end{itemize}}
\end{rmk}

We now deduce that for compactly aligned $(\Lambda,d)$, the groupoid $G_\Lambda$ has a locally compact topology; recall from Remark~\ref{rmk:not compactly aligned entails problems} that if $(\Lambda,d)$ fails to be compactly aligned, then the topology on $G_\Lambda$ may not be locally compact.

\begin{prop}\label{prop:topology on G_Lambda is locally compact}
Let $(\Lambda,d)$ be a compactly aligned topological $k$-graph. For $p,q\in\NN^k$ and compact sets $U\subset\Lambda^p$ and $V\subset\Lambda^q$, the set $Z(U*_sV,p-q)$ is compact.
\end{prop}

\begin{proof}
Let $\langle (\lambda_jx_j,p-q,\mu_jx_j)\rangle_{j\in\NN}$ be a sequence in $Z(U*_sV,p-q)$ with each $(\lambda_j,\mu_j)\in U*_sV$. Since $U*_sV$ is compact, there exists an infinite set $I_0\subset\NN$ such that $\langle (\lambda_j,\mu_j)\rangle_{j\in I_0}$ converges to some $(\lambda,\mu)\in U*_sV$. We construct an element $z\in X_\Lambda$ and an infinite set $I\subset I_0$ such that $\langle (\lambda_jx_j,p-q,\mu_jx_j)\rangle_{j\in I}$ converges to $(\lambda z,p-q,\mu z)$.

Define $f:\NN\to\{1,\dots,k\}$ by $f(j)=j (\modulo k)+1$, and iteratively construct a sequence $\langle z_j\rangle_{j\in\NN}$ in $\Lambda$ as follows: First, set $z_0:=s(\lambda)=s(\mu)$. Let $l\in\NN$, and suppose $z_0,\dots,z_l\in\Lambda$ and infinite sets $I_l\subset\cdots\subset I_0\subset\NN$ have been defined and satisfy:
\begin{equation}
s(z_i)=r(z_{i+1}) \quad\text{for all } 0\le i\le l-1,
\end{equation}
\begin{equation}
d(x_j)\ge d(z_0\cdots z_l) \quad\text{for all } j\in I_l,
\end{equation}
and
\begin{equation}
\lim_{j\in I_l}(\lambda_jx_j(0,d(z_0\cdots z_l)),\mu_jx_j(0,d(z_0\cdots z_l)))=(\lambda z_0\cdots z_l, \mu z_0\cdots z_l).
\end{equation}
One of the following two properties must hold:
\begin{itemize}
\item[1)]
There exists a compact set $W_{l+1}\subset\Lambda^{e_{f(l+1)}}$ and an infinite set $I'_{l+1}\subset I_l$  such that $d(x_j)\ge d(z_0\cdots z_l) + e_{f(l+1)}$ for all $j\in I'_{l+1}$, and the sequence $\langle x_j(d(z_0\cdots z_l),d(z_0\cdots z_l) + e_{f(l+1)})\rangle_{j\in I'_{l+1}}$ is contained in $W_{l+1}$;
\item[2)]
The set $J_{l+1}:=\{j\in I_l ~|~ d(x_j)\ge d(z_0\cdots z_l) + e_{f(l+1)}\}$ is either
\begin{itemize}
\item[-]
finite, or
\item[-]
infinite and the sequence $\langle x_j(d(z_0\cdots z_l),d(z_0\cdots z_l) + e_{f(l+1)})\rangle_{j\in J_{l+1}}$ is wandering.
\end{itemize}
\end{itemize}
If 1) holds, then fix compact $W_{l+1}$ and infinite $I'_{l+1}$ satisfying the conditions of 1); as $W_{l+1}$ is compact, fix an infinite subset $I_{l+1}\subset I'_{l+1}$ such that the sequence $\langle x_j(d(z_0\cdots z_l),d(z_0\cdots z_l) + e_{f(l+1)})\rangle_{j\in I_{l+1}}$ converges in $W_{l+1}$, and define
\[
z_{l+1}:=\lim_{j\in I_{l+1}}x_j(d(z_0\cdots z_l),d(z_0\cdots z_l) + e_{f(l+1)}).
\]
If 2) holds, then set $z_{l+1}:=s(z_l)$ and $I_{l+1}:=I_l$.

There exists a unique path $z\in X_\Lambda$ such that
\[
d(z) = \lim_{l\to\infty}d(z_0\cdots z_l) \quad\text{and}\quad z(0,d(z_0\cdots z_l)) = z_0\cdots z_l \text{ for all } l\in\NN.
\]
We also have for all $l\in\NN$,
\begin{equation}\label{eqn:initial segments converge}
\lim_{j\in I_l} (\lambda_j x_j(0,d(z_0\cdots z_l)),\mu_jx_j(0,d(z_0\cdots z_l))) = (\lambda z_0\cdots z_l,\mu z_0\cdots z_l)
\end{equation}
and
\begin{equation}\label{eqn:property of z}
\begin{split}
&\text{if } J_{l+1}:=\{j\in I_l ~|~ d(x_j)\ge d(z_0\cdots z_l) + e_{f(l+1)}\} \text{ is infinite}\\
&\text{and } \langle x_j(d(z_0\cdots z_l),d(z_0\cdots z_l)+e_{f(l+1)})\rangle_{j\in J_{l+1}} \text{ is not wandering}, \\
&\text{then } z_{l+1}\in\Lambda^{e_{f(l+1)}}; \text{ otherwise } z_{l+1}=s(z_l).
\end{split}
\end{equation}

Define an infinite set $I = \{j_i\}_{i\in\NN}\subset\NN$ by choosing any $j_0\in I_0$, and, after $j_0,\dots,j_l$ have been set, choosing $j_{l+1}\in I_{l+1}$ such that $j_{l+1}\ge j_l$.

We claim that Conditions~(1) and~(2) of Proposition~\ref{prop:convergence means} hold for the sequence $\langle (\lambda_jx_j,p-q,\mu_jx_j)\rangle_{j\in I}$ and point $(\lambda z,p-q,\mu z)$.

Condition~(1) of Proposition~\ref{prop:convergence means} follows from \eqref{eqn:initial segments converge}. Suppose for contradiction that Condition~(2) of Proposition~\ref{prop:convergence means} does not hold, so there exist $i\in\{1,\dots,k\}$ and $n\in\NN^k$ such that $n\le d(z)$, $n_i=d(z)_i$,
\[
J(p+n,i):=\{j\in I ~|~ d(\lambda_jx_j)\ge p+n+e_i\}
\]
is infinite, and the sequence $\langle (\lambda_jx_j)(p+n,p+n+e_i)\rangle_{j\in J(p+n,i)}$ is not wandering. Then there exists a compact set $V_1\subset\Lambda^{e_i}$ such that $V_1$ contains infinitely many elements of $\langle (\lambda_jx_j)(p+n,p+n+e_i)\rangle_{j\in J(p+n,i)}$.

Let $L\in\NN$ be the smallest number such that $d(z_0\cdots z_L)\ge n$ and $f(L+1)=i$, so $d(z)_i=n_i=d(z_0\cdots z_L)_i$. Let $V_2\subset\Lambda^{d(z_0\cdots z_L)-n}$ be a compact neighbourhood of $z(n,d(z_0\cdots z_L))$. Then $\langle x_j(n,d(z_0\cdots z_L))\rangle_{j\in J(p+n,i)}$ is eventually in $V_2$ by Condition~(1).

Since $\Lambda$ is compactly aligned and $(d(z_0\cdots z_L)-n)_i=0$, it follows that $V_1\vee V_2$ is a compact subset of $\Lambda^{d(z_0\cdots z_L)-n+e_i}$. Furthermore, $V_1\vee V_2$ contains infinitely many elements of $\langle x_j(n,d(z_0\cdots z_L)+e_i)\rangle_{j\in J(p+n,i)}$. Since $J(p+n,i)$ is contained in the set $J_{L+1}$ of \eqref{eqn:property of z}, we have that $\Seg^{d(z_0\cdots z_L)-n+e_i}_{(d(z_0\cdots z_L)-n, d(z_0\cdots z_L)-n+e_i)}(V_1\vee V_2)$ is a compact subset of $\Lambda^{e_i}$ which contains infinitely many elements of $\langle x_j(d(z_0\cdots z_L),d(z_0\cdots z_L)+e_i)\rangle_{j\in J_{L+1}}$. By \eqref{eqn:property of z}, we then have $z_{L+1}\in\Lambda^{e_i}$, which contradicts $d(z_0\cdots z_L)_i=d(z)_i$. Therefore Condition~(2) of Proposition~\ref{prop:convergence means} holds, and Proposition~\ref{prop:convergence means} implies
\[
\lim_{j\in I} (\lambda_jx_j,p-q,\mu_jx_j) = (\lambda z,p-q,\mu z),
\]
completing the proof.
\end{proof}

\begin{theorem}
Let $(\Lambda,d)$ be a compactly aligned topological $k$-graph. Then $G_\Lambda$ is a locally compact $r$-discrete topological groupoid admitting a Haar system consisting of counting measures.
\end{theorem}

\begin{proof}
Local compactness of $G_\Lambda$ follows from Proposition~\ref{prop:topology on G_Lambda is locally compact}. Straightforward applications of Proposition~\ref{prop:convergence means} give continuity of composition and inversion. Therefore $G_\Lambda$ is a locally compact topological groupoid.

To show that $G_\Lambda$ is $r$-discrete and admits a Haar system, by \cite[Proposition~I.2.8]{Ren1} it suffices to show that $r:G_\Lambda\to G_\Lambda^{(0)}$ is a local homeomorphism. Fixing $(\lambda x,d(\lambda)-d(\mu),\mu x)\in G_\Lambda$ and choosing open neighbourhoods $U\subset\Lambda^{d(\lambda)}$ and $V\subset\Lambda^{d(\mu)}$ of $\lambda$ and $\mu$, respectively, such that $s|_U$ and $s|_V$ are homeomorphisms, one checks that $r|_{Z(U*_sV,d(\lambda)-d(\mu))}$ is a homeomorphism. Therefore $G_\Lambda$ is $r$-discrete and admits a Haar system, and by \cite[Lemma~I.2.7]{Ren1} we can choose the Haar system to comprise counting measures.
\end{proof}

\begin{examples}\label{ex:path-groupoid examples}
\begin{enumerate}
\item\label{ex:item:directed graph path groupoid}
Let $E$ be a discrete directed graph, and recall the construction of the topological $1$-graph $(E^*,l)$ from Example~\ref{ex:graph examples}\eqref{ex:item:top graph}. In \cite{Pat2}, Paterson defines an inverse semigroup $S^{\rm Pat}_E$ and an action of $S^{\rm Pat}_E$ on the path space $X_{E^*}$. He then defines the topological groupoid $H^{\rm Pat}_E$ as the groupoid of germs of the action. Comparing $G_{E^*}$ with the description of $H^{\rm Pat}_E$ given in \cite[Theorem~1]{Pat2} and comparing the topological structures of $G_{E^*}$ and $H^{\rm Pat}_E$ given in Proposition~\ref{prop:defining a topology on G_Lambda} and \cite[Proposition~3]{Pat2}, respectively, we see that $G_{E^*}$ and $H^{\rm Pat}_E$ are isomorphic as topological groupoids.
\item\label{ex:item:k-graph path groupoid}
Given a finitely aligned discrete $k$-graph, the authors of \cite{FMY1} define and study an $r$-discrete groupoid $G^{\rm FMY}_\Lambda$ (see \cite[Section~6]{FMY1}). Comparing $G_\Lambda$ with the description of $G^{\rm FMY}_\Lambda$ given in \cite[Remark~6.2]{FMY1} and comparing the topological structures of $G_\Lambda$ and $G^{\rm FMY}_\Lambda$ given in Proposition~\ref{prop:defining a topology on G_Lambda} and \cite[Remark~6.4]{FMY1}, respectively, we see that $G_\Lambda$ and $G^{\rm FMY}_\Lambda$ are isomorphic as topological groupoids.
\end{enumerate}
\end{examples}

\section{The boundary-path groupoid}\label{sec:boundary-path gpd}

Given a compactly aligned topological $k$-graph $(\Lambda,d)$, we now identify a closed invariant subset $\partial\Lambda$ of $X_\Lambda = G_\Lambda^{(0)}$, and define our boundary-path groupoid $\Gg_\Lambda$ as the reduction $G_\Lambda|_{\partial\Lambda}$.

\begin{dfn}
Let $(\Lambda,d)$ be a topological $k$-graph and let $V\subset\Lambda^0$. A set $E\subset V\Lambda$ is \emph{exhaustive} for $V$ if for all $\lambda\in V\Lambda$ there exists $\mu\in E$ such that $\Lmin(\lambda,\mu)\neq\emptyset$. For $v\in\Lambda^0$, let $v\CE(\Lambda)$ denote the set of all compact sets $E\subset\Lambda$ such that $r(E)$ is a neighbourhood of $v$ and $E$ is exhaustive for $r(E)$.
\end{dfn}

\begin{dfn}
Let $(\Lambda,d)$ be a topological $k$-graph. A path $x\in X_\Lambda$ is called a \emph{boundary path} if for all $m\in\NN^k$ with $m\le d(x)$, and for all $E\in x(m)\CE(\Lambda)$, there exists $\lambda\in E$ such that $x(m,m+d(\lambda))=\lambda$. We write $\partial\Lambda$ for the set of all boundary paths in $X_\Lambda$. For $v\in\Lambda^0$ and $V\subset\Lambda^0$, we define $v\partial\Lambda = \{x\in\partial\Lambda ~|~ r(x) = v\}$ and $V(\partial\Lambda) = \{x\in\partial\Lambda ~|~ r(x)\in V\}$.
\end{dfn}

\begin{prop}\label{prop:partialLambda is nonempty}
Let $(\Lambda,d)$ be a topological $k$-graph. Then $v\partial\Lambda$ is nonempty for all $v\in\Lambda^0$.
\end{prop}

\begin{proof}
We construct a path $x\in v\partial\Lambda$. Define $f:\NN\to\{1,\dots,k\}$ by $f(j):= j (\modulo k)+1$. If $v\Lambda^{e_{f(1)}}\neq \emptyset$, then choose $\lambda_1\in v\Lambda^{e_{f(1)}}$, otherwise set $\lambda_1:=v$. Once $\lambda_1,\dots,\lambda_l$ have been defined, choose $\lambda_{l+1}\in s(\lambda_l)\Lambda^{e_{f(l+1)}}$ if $s(\lambda_l)\Lambda^{e_{f(l+1)}}\neq\emptyset$, otherwise set $\lambda_{l+1}:= s(\lambda_l)$. There exists $x\in X_\Lambda$ such that $d(x) =\lim_{j\in\NN} d(\lambda_1\cdots\lambda_j)$ and $x(0,d(\lambda_1\cdots\lambda_j))=\lambda_1\cdots\lambda_j$ for all $j\in\NN$.

To show $x\in\partial\Lambda$, let $m\in\NN^k$ satisfy $m\le d(x)$ and let $E\in x(m)\CE(\Lambda)$. Since $E$ is exhaustive, for each $n\in\NN^k$ with $m\le n\le d(x)$, there exists $\mu_n\in E$ such that $\Lmin(x(m,n),\mu_n)\neq\emptyset$. We will show there exists $N\in\NN^k$ such that $d(\mu_N)\le N-m$; for this $N$ we then have $x(m,m+d(\mu_N))=\mu_N$, as required.

Since $E$ is compact, it follows that $\{d(\mu_n) ~|~ m\le n\le d(x)\}$ is finite. If $d(x)_i=\infty$ for all $i\in\{1,\dots,k\}$, then choosing $N:= \bigvee \{m+d(\mu_n) ~|~ m\le n\le d(x)\}$ will do. So suppose there exists $i\in\{1,\dots,k\}$ such that $d(x)_i<\infty$.

Define $I:=\{i\in\{1,\dots,k\} ~|~ d(x)_i<\infty\}$ and
\begin{equation}\label{eqn:defining p}
p:= \big(\textstyle{\bigvee}\{m+d(\mu_n) ~|~ m\le n\le d(x)\}\big)\wedge d(x).
\end{equation}
For each $i\in I$, let $l_i\in\NN$ be the smallest number such that $f(l_i)=i$ and $s(\lambda_{i-1})\Lambda^{e_i}=\emptyset$, so $s(\lambda_j)\Lambda^{e_i}=\emptyset$ for all $j\ge l_i$. Let $L\in\NN$ be the smallest number such that
\begin{equation}\label{eqn:defining L}
L\ge\max_{i\in I}l_i \quad\text{and}\quad d(\lambda_1\cdots\lambda_L)\ge p.
\end{equation}

Define $N:=d(\lambda_1\cdots\lambda_L)$. Suppose for contradiction that $d(\mu_N)\not\le N-m$, so there exists $i\in\{1,\dots,k\}$ such that $d(\mu_N)_i>(N-m)_i$. Then~\eqref{eqn:defining p} and~\eqref{eqn:defining L} imply $N_i=d(x)_i$. Thus $d(x)_i<\infty$, and it follows from \eqref{eqn:defining L} that $s(\lambda_L)\Lambda^{e_i}=\emptyset$. But $\Lmin(x(m,N),\mu_N)\neq\emptyset$ and for any $(\alpha,\beta)\in\Lmin(x(m,N),\mu_N)$ we must have $d(\alpha)_i > 0$, contradicting $s(\lambda_L)\Lambda^{e_i}=\emptyset$. Therefore $d(\mu_N)\le N-m$, so $x(m,m+d(\mu_N))=\mu_N$ and $x\in v\partial\Lambda$.
\end{proof}

\begin{prop}\label{prop:partialLambda is closed}
Let $(\Lambda,d)$ be a topological $k$-graph. $\partial\Lambda$ is closed in $G_\Lambda^{(0)}$.
\end{prop}

\begin{proof}
Let $\langle x_j\rangle_{j\in\NN}$ be a sequence in $\partial\Lambda$ converging to some $x\in X_\Lambda$. Suppose for contradiction that $x\not\in\partial\Lambda$, so there exists $m\in\NN^k$, $m\le d(x)$, and $E\in x(m)\CE(\Lambda)$ such that $x(m,p)\not\in E$ for all $p\in\NN^k$ with $m\le p\le d(x)$.

Let $U\subset \Lambda^m$ be a relatively compact open neighbourhood of $x(0,m)$ such that $s(U)\subset r(E)$. Then $x\in Z(U)\cap Z(\overline{U}E)^c$, so there exists $J\in\NN$ such that $x_j\in Z(U)\cap Z(\overline{U}E)^c$ whenever $j\ge J$. But then for $j\ge J$ and $p\in\NN^k$ with $m\le p\le d(x_j)$, we have $x_j(0,p)\not\in \overline{U}E$, which implies $x_j(m,p)\not\in E$, contradicting $x_j\in\partial\Lambda$ and $E\in x_j(m)\CE(\Lambda)$. Hence $x\in\partial\Lambda$, and $\partial\Lambda$ is closed.
\end{proof}

To prove $\partial\Lambda$ is an invariant subset of $G_\Lambda^{(0)}$, we first need a definition and a lemma.

\begin{dfn}
Let $(\Lambda,d)$ be a topological $k$-graph. For $E,F\subset\Lambda$, define the \emph{minimal extenders of $E$ by $F$} to be the set
\[
\Ext(E;F) = \bigcup_{\lambda\in E}\bigcup_{\mu\in F}\{\alpha\in\Lambda ~|~ (\alpha,\beta)\in \Lmin(\lambda,\mu) \text{ for some } \beta\in\Lambda\}.
\]
If $E$ is a singleton set $E=\{\lambda\}$, we write $\Ext(\lambda;F)$ for $\Ext(\{\lambda\};F)$.
\end{dfn}

The proof of the following lemma is straightforward.

\begin{lemma}[{cf. \cite[Lemma~C.5]{RSY2}}]\label{lem:CE sets are transitive}
Let $(\Lambda,d)$ be a compactly aligned topological $k$-graph, let $v\in\Lambda^0$ and $\lambda\in v\Lambda$, and suppose $E\in v\CE(\Lambda)$. Then for any compact neighbourhood $U\subset\Lambda^{d(\lambda)}$ of $\lambda$, $\Ext(U;E)\in s(\lambda)\CE(\Lambda)$.
\end{lemma}

\begin{prop}\label{prop:partialLambda invariant}
Let $(\Lambda,d)$ be a compactly aligned topological $k$-graph. For $x\in\partial\Lambda$, $m\in\NN^k$ with $m\le d(x)$, and $\lambda\in \Lambda r(x)$, we have $\sigma^mx, \lambda x\in\partial\Lambda$. Hence $\partial\Lambda$ is an invariant subset of $G_\Lambda^{(0)}$.
\end{prop}

\begin{proof}
To see that $\sigma^mx\in\partial\Lambda$, let $n\in\NN^k$ satisfy $n\le d(\sigma^mx)$, and let $E\in (\sigma^mx)(n)\CE(\Lambda)$. Then $E\in x(m+n)\CE(\Lambda)$, so there exists $\mu\in E$ such that $(\sigma^mx)(n,n+d(\mu))=x(m+n,m+n+d(\mu))=\mu$, as required.

Now let $n\in\NN^k$ satisfy $n\le d(\lambda x)$, and let $E\in (\lambda x)(n)\CE(\Lambda)$. Define $\xi := x(n,n\vee d(\lambda))$ and let $U\subset \Lambda^{d(\xi)}$ be a relatively compact open neighbourhood of $\xi$ such that $s|_U$ is a homeomorphism. Then by Lemma~\ref{lem:CE sets are transitive} we have
\[
\Ext(\overline{U};E)\in (\lambda x)(n\vee d(\lambda))\CE(\Lambda) = x((n\vee d(\lambda))-d(\lambda))\CE(\Lambda),
\]
so there exists $\alpha\in\Ext(\overline{U};E)$ such that
\[
x((n\vee d(\lambda))-d(\lambda),(n\vee d(\lambda))-d(\lambda)+d(\alpha)) = \alpha.
\]
Since $s(\xi)=r(\alpha)$ and $s|_{U}$ is a homeomorphism, we have $\xi\alpha=\mu\beta$ for some $\mu\in E$ and $\beta\in\Lambda$. Hence $(\lambda x)(n,n+d(\mu))=\mu$, giving $\lambda x\in\partial\Lambda$.
\end{proof}

\begin{dfn}
Let $(\Lambda,d)$ be a compactly aligned topological $k$-graph. Propositions~\ref{prop:partialLambda is nonempty},~\ref{prop:partialLambda is closed} and~\ref{prop:partialLambda invariant} imply that $\partial\Lambda$ is a nonempty closed invariant subset of $G_\Lambda^{(0)}$, and we define the \emph{boundary-path groupoid} $\Gg_\Lambda$ to be the reduction $\Gg_\Lambda := G_\Lambda|_{\partial\Lambda}$. Then $\Gg_\Lambda$ is a locally compact $r$-discrete groupoid admitting a Haar system consisting of counting measures.
\end{dfn}

\begin{notation}\label{notation:basis sets for bp space}
To distinguish basis sets of $\Gg_\Lambda$ from those of $G_\Lambda$, for $F\subset\Lambda*_s\Lambda$, $m\in\ZZ^k$ and $U\subset\Lambda$, define
\[
\Zz(F,m) = Z(F,m)\cap \Gg_\Lambda
\]
and
\[
\Zz(U) = Z(U)\cap \Gg_\Lambda = \Zz(U*_sU,0)\cap\Zz(\Lambda^0*_s\Lambda^0,0).
\]
\end{notation}

\begin{examples}\label{ex:boundary-path groupoids}
\begin{enumerate}
\item\label{ex:item:directed graph boundary-path groupoid}
Recalling Example~\ref{ex:path-groupoid examples}\eqref{ex:item:directed graph path groupoid}, for a discrete directed graph $E$, Paterson \cite{Pat2} defines an $r$-discrete groupoid $H^{\rm Pat}_E$, and we saw that $G_{E^*}$ and $H^{\rm Pat}_E$ are isomorphic. Paterson then identifies a closed invariant subset $X^{\rm Pat}$ of $(H^{\rm Pat}_E)^{(0)}$ (see paragraph preceding \cite[Proposition~5]{Pat2}) and studies the reduction of $H^{\rm Pat}_E$ by $X^{\rm Pat}$ (see \cite[Theorem~2]{Pat2}). In the setting of directed graphs, it is straightforward to see that our boundary paths are precisely the elements of Paterson's set $X^{\rm Pat}$, so it follows that $\Gg_{E^*}$ and $H^{\rm Pat}_E|_{X^{\rm Pat}}$ are isomorphic as topological groupoids.
\item\label{ex:item:k-graph boundary-path groupoid}
Let $(\Lambda,d)$ be a finitely aligned discrete $k$-graph. As discussed in Example~\ref{ex:path-groupoid examples}\eqref{ex:item:k-graph path groupoid}, the authors of \cite{FMY1} associate a topological groupoid $G^{\rm FMY}_\Lambda$ to $\Lambda$, and this groupoid is isomorphic to our path groupoid $G_\Lambda$. In \cite{FMY1}, the authors identify a set of boundary paths $\partial\Lambda$ of $\Lambda$ as a closed invariant subset of $(G^{\rm FMY}_\Lambda)^{(0)}$ and study the reduction $G^{\rm FMY}_\Lambda|_{\partial\Lambda}$. Since our definition of a boundary path corresponds to \cite[Definition~5.10]{FMY1} in the setting of finitely aligned discrete $k$-graphs, it follows that $\Gg_\Lambda$ and $G^{\rm FMY}_\Lambda|_{\partial\Lambda}$ are isomorphic as topological groupoids.
\item\label{ex:item:SGDS groupoid}
Given a singly generated dynamical system $(X,\theta)$, we formed a topological $1$-graph $(\Lambda(X,\theta),d)$ (see Example~\ref{ex:graph examples}\eqref{ex:item:SGDS}). In \cite{Ren4}, Renault defines a topological groupoid $G(X,\theta)\subset X\times\ZZ\times X$ by setting
\[
G(X,\theta) = \{(x,m-n,y) ~|~ x\in\dom(\theta^m), y\in\dom(\theta^n), \theta^m(x) = \theta^n(y)\},
\]
with the usual groupoid structure, and defining basis sets
\[
\Uu(U;m,n;V) = \{(x,m-n,y) ~|~ (x,y)\in U\times V, \theta^m(x)=\theta^n(y)\}
\]
where $U\subset\dom(\theta^m)$, $V\subset\dom(\theta^n)$ are open sets on which, respectively, $\theta^m$ and $\theta^n$ are injective (see \cite[Section~2]{Ren4} for details; see also \cite{De1} for the same construction with $X$ compact and $\theta$ surjective).

The boundary paths of $\Lambda(X,\theta)$ can be identified with $X$, and under this identification, $\sigma$ is intertwined with $\theta$. One can then show that the groupoid and topological structures on $\Gg_{\Lambda(X,\theta)}$ and $G(X,\theta)$ are equivalent, hence the two topological groupoids are isomorphic.
\end{enumerate}
\end{examples}

\section[Aperiodicity and essential freeness]{Aperiodicity in topological higher-rank graphs and essential freeness of boundary-path groupoids}\label{sec:aperiodicity}

In this section we consider an analogue of the Aperiodicity Condition used in \cite{KP2,FMY1}. Using the condition, we extend \cite[Proposition~4.5]{KP2} and \cite[Proposition~7.2]{FMY1} to our setting.

\begin{dfn}
Let $(\Lambda,d)$ be a topological $k$-graph. A boundary path $x\in\partial\Lambda$ is \emph{aperiodic} if
\begin{equation}\label{eqn:aperiodic path}
\text{for all } p,q\in\NN^k \text{ with } p,q\le d(x),~ p\neq q \text{ implies } \sigma^px\neq \sigma^q x.
\end{equation}
\end{dfn}

Recall that a topological groupoid $\Gamma$ is \emph{essentially free} if the set of units with trivial isotropy is dense in $\Gamma^{(0)}$; that is, $\overline{\{ x\in\Gamma^{(0)} ~|~ x\Gamma x = \{x\}\}}=\Gamma^{(0)}$.

\begin{theorem}\label{thm:aperiodicity iff essentially free}
Let $(\Lambda,d)$ be a compactly aligned topological $k$-graph. Then $\Gg_\Lambda$ is essentially free if and only if
\begin{equation}\label{eqn:aperiodicity condition}\tag{A}
\text{for all nonempty open } V\subset\Lambda^0,\text{ there exists an aperiodic path } x\in V(\partial\Lambda).
\end{equation}
\end{theorem}

To prove the theorem, we need the following two lemmas.

\begin{lemma}\label{lem:aperiodic iff trivial isotropy}
Let $(\Lambda,d)$ be a compactly aligned topological $k$-graph. A boundary-path $x\in\partial\Lambda$ is aperiodic if and only if its associated isotropy group in $\Gg_\Lambda$ is trivial.
\end{lemma}

\begin{proof}
The lemma follows from the equivalence: For $m\in\ZZ^k$, the triple $(x,m,x)$ is an element of $\Gg_\Lambda$ if and only if there exist $p,q\in\NN^k$ such that $p,q\le d(x)$, $p-q=m$ and $\sigma^p x=\sigma^q x$.
\end{proof}

\begin{lemma}\label{lem:aperiodicity is transitive}
Let $(\Lambda,d)$ be a topological $k$-graph. For any aperiodic $x\in\partial\Lambda$ and $\lambda\in\Lambda r(x)$, $\lambda x$ is aperiodic.
\end{lemma}

\begin{proof}
Arguing by contrapositive, suppose that $\lambda x$ is not aperiodic, so there exists $p,q\in\NN^k$ such that $p,q\le d(\lambda x)$, $p\neq q$ and $\sigma^p(\lambda x)\neq \sigma^q(\lambda x)$. It follows that
\[
(d(\lambda)+p)\wedge d(\lambda x)\neq (d(\lambda) +q)\wedge d(\lambda x)
\]
and
\[
\sigma^{(d(\lambda)+p)\wedge d(\lambda x)}(\lambda x)=\sigma^{(d(\lambda) +q)\wedge d(\lambda x)}(\lambda x).
\]
Thus we have $p\wedge d(x)\neq q\wedge d(x)$ and
\[
\sigma^{p\wedge d(x)}x = \sigma^{(d(\lambda)+p)\wedge d(\lambda x)}(\lambda x)=\sigma^{(d(\lambda) +q)\wedge d(\lambda x)}(\lambda x) = \sigma^{q\wedge d(x)}x,
\]
proving $x$ is not aperiodic.
\end{proof}

\begin{proof}[{Proof of Theorem~\ref{thm:aperiodicity iff essentially free}}]

First assume that $\Gg_\Lambda$ is essentially free, and let $V\subset\Lambda^0$ be nonempty and open. By Lemma~\ref{prop:partialLambda is nonempty}, $V(\partial\Lambda)$ is nonempty, so the open set $Z(V)$ is nonempty in $\Gg_\Lambda^{(0)}$. Therefore, there exists $x\in Z(V)$ with trivial isotropy, and Lemma~\ref{lem:aperiodic iff trivial isotropy} implies that $x\in V(\partial\Lambda)$ is aperiodic. Hence $(\Lambda,d)$ satisfies Condition~\eqref{eqn:aperiodicity condition}.

Conversely, assume that Condition~\eqref{eqn:aperiodicity condition} holds. Fix $x\in \Gg_\Lambda^{(0)}$ and let $\Zz(U)\cap \Zz(F)^c$ be a basis set containing $x$. There exists $\lambda\in U$ such that $x(0,d(\lambda))=\lambda$; we can assume $U\subset\Lambda^{d(\lambda)}$, $U$ is relatively compact and open, $s|_U$ is a homeomorphism, and, by Lemma~\ref{lem:choosing better basis sets for unit space}, every $\mu\in F$ has the form $\mu=\xi\mu'$ for some $\xi\in\overline{U}$.

The set $\{d(\mu) ~|~ \mu\in F\}$ is finite, and $d(\lambda)\le d(\mu)$ for all $\mu\in F$, so we define the compact set
\[
E := \bigcup_{m\in\{d(\mu) ~|~ \mu\in F\}}\Seg^m_{(d(\lambda),m)}(F\cap \Lambda^m).
\]
We know that $(s|_{\overline{U}})^{-1}(s(\lambda))=\{\lambda\}$, so for $\nu\in E$, if $r(\nu)=s(\lambda)$, then $\lambda\nu\in F$. It follows that if there exists $\nu\in E$ such that $x(d(\lambda),d(\lambda)+d(\nu))=\nu$, then $x(0,d(\lambda\nu))\in F$, which contradicts $x\in \Zz(F)^c$. Since $x$ is a boundary path, we deduce that $E\not\in s(\lambda)\CE(\Lambda)$.

Now, $E$ must fail to be an element of $s(\lambda)\CE(\Lambda)$ on account of one of two reasons: either $r(E)$ is not a neighbourhood of $s(\lambda)$, or $E$ is not exhaustive for $r(E)$. In either case, there exists $\eta\in s(U)\Lambda$ such that $\Lmin(\eta,\nu)=\emptyset$ for all $\nu\in E$.

We claim there exists a neighbourhood $W\subset\Lambda^{d(\eta)}$ of $\eta$ such that $\Lmin(\eta',\nu)=\emptyset$ for all $\eta'\in W$ and $\nu\in E$: Suppose for contradiction that there exist sequences $\langle\eta_j\rangle_{j\in\NN}\subset\Lambda^{d(\eta)}$ and $\langle\nu_j\rangle_{j\in\NN}\subset E$ such that $\lim_{j\in\NN}\eta_j=\eta$ and $\Lmin(\eta_j,\nu_j)\neq\emptyset$ for each $j\in\NN$. Since $E$ is compact, we can assume $\langle\nu_j\rangle_{j\in\NN}$ converges to some $\nu\in E$.

Let $\langle(\alpha_j,\beta_j)\rangle_{j\in\NN}$ be a sequence such that $(\alpha_j,\beta_j)\in\Lmin(\eta_j,\nu_j)$ for each $j\in\NN$. We can assume $d(\nu_j)=d(\nu)$ for all $j\in\NN$, and hence $d(\alpha_j)=(d(\eta)\vee d(\nu))-d(\eta)$ and $d(\beta_j)=(d(\eta)\vee d(\nu))-d(\nu)$ for all $j\in\NN$.

Taking a compact neighbourhood $Y\subset\Lambda^{d(\eta)}$ of $\eta$, we can assume $\langle\eta_j\rangle_{j\in\NN}\subset Y$. Since $\Lambda$ is compactly aligned, it follows that $Y\vee E$ is compact and contains the sequence $\langle\eta_j\alpha_j\rangle_{j\in\NN} = \langle\nu_j\beta_j\rangle_{j\in\NN}$.
Hence there exists a convergent subsequence $\langle\eta_j\alpha_j\rangle_{j\in J}$ contained in $\Lambda^{d(\eta)\vee d(\nu)}$. We must then have $\lim_{j\in J}\eta_j\alpha_j=\eta\alpha$ and $\lim_{j\in J}\nu_j\beta_j=\nu\beta$ for some $\alpha,\beta\in\Lambda$.
It follows that $(\alpha,\beta)\in\Lmin(\eta,\nu)$; a contradiction. Therefore there exists a neighbourhood $W\subset\Lambda^{d(\eta)}$ of $\eta$ such that $\Lmin(\eta',\nu)=\emptyset$ for all $\eta'\in W$ and $\nu\in E$. We can assume that $r(W)\subset s(U)$.

By Condition~(A), there exists an aperiodic path $z\in s(W)(\partial\Lambda)$. Since $r(W)\subset s(U)$, there exist $\lambda'\in U$ and $\eta'\in W$ such that $s(\lambda')=r(\eta')$ and $s(\eta')=r(z)$. Proposition~\ref{prop:partialLambda invariant} and Lemma~\ref{lem:aperiodicity is transitive} imply $\lambda' \eta' z$ is an aperiodic boundary path. Since $\Lmin(\eta',\nu)=\emptyset$ for all $\nu\in E$, it follows that $\lambda'\eta' z\in \Zz(U)\cap \Zz(F)^c$. Lemma~\ref{lem:aperiodic iff trivial isotropy} now gives the result.
\end{proof}

\section[{Amenability of $\Gg_\Lambda$}]{Amenability of the boundary-path groupoid}\label{sec:amenability}

In this section we prove amenability of the boundary-path groupoid under certain conditions on the topological $k$-graph. Rather than detail the characterizations of groupoid amenability here, we refer the reader to \cite[Chapter~2]{A-DRen}.

Recall that a locally compact groupoid $\Gamma$ is \emph{proper} if $(r,s):\Gamma\to\Gamma^{(0)}\times\Gamma^{(0)}$, defined by $\gamma\mapsto (r(\gamma),s(\gamma))$, is a proper mapping; that is, if the inverse image of any compact set from $\Gamma^{(0)}\times\Gamma^{(0)}$ is compact.

The next proposition is standard in groupoid theory; we omit its proof.

\begin{prop}\label{prop:proper implies amenable}
Let $\Gamma$ be a locally compact proper groupoid admitting a Haar system. Then $\Gamma$ is amenable in the sense of \cite[Definition~2.2.7]{A-DRen}.
\end{prop}

\begin{prop}
Let $(\Lambda,d)$ be a topological $k$-graph. If either
\begin{itemize}
\item[(i)]
$k=1$, or
\item[(ii)]
$(\Lambda,d)$ is a finitely aligned discrete $k$-graph,
\end{itemize}
then $\Gg_\Lambda$ is amenable.
\end{prop}

\begin{proof}
Suppose $k=1$. Then \cite[Theorem~5.2]{Y2} implies $C^*(\Gg_\Lambda)$ is isomorphic to the Cuntz-Pimsner algebra $\Oo(E_\Lambda)$ (see \cite[Section~5]{Y2}). Thus, \cite[Proposition~6.1]{Kat1} implies $C^*(\Gg_\Lambda)$ is nuclear, and since $\Gg_\Lambda$ is $r$-discrete, \cite[Corollary~6.2.14 and Theorem~3.3.7]{A-DRen} imply $\Gg_\Lambda$ is amenable.

Now suppose $(\Lambda,d)$ is a finitely aligned discrete $k$-graph. By \cite[Theorem~6.13]{FMY1}, $C^*(\Gg_\Lambda)\cong C^*(\Lambda)$, where $C^*(\Lambda)$ is defined in \cite[Remark~3.9]{FMY1}. We then deduce $C^*(\Gg_\Lambda)$ is nuclear by \cite[Proposition~8.1]{S2}. Since $\Gg_\Lambda$ is $r$-discrete, \cite[Corollary~6.2.14 and Theorem~3.3.7]{A-DRen} imply $\Gg_\Lambda$ is amenable.
\end{proof}

\begin{dfn}
Let $(\Lambda,d)$ be a topological $k$-graph and let $v\in\Lambda^0$. Then $v$ is said to be a \emph{source} if $v\Lambda^{e_i}=\emptyset$ for some $i\in\{1,\dots,k\}$, and $v$ is said to be a sink if $\Lambda^{e_i}v=\emptyset$ for some $i\in\{1,\dots,k\}$.
\end{dfn}

The following definition generalizes the row-finite condition on discrete higher-rank graphs.

\begin{dfn}
A topological $k$-graph is \emph{proper} if, for all $m\in\NN^k$, $r|_{\Lambda^m}$ is a proper map; that is, if, for all $m\in\NN^k$ and compact $U\subset\Lambda^0$, $U\Lambda^m$ is compact.
\end{dfn}

\begin{rmk}
It is straightforward to see that any proper topological $k$-graph is compactly aligned.
\end{rmk}

\begin{lemma}
Let $(\Lambda,d)$ be a proper topological $k$-graph without sources. Then $d(x) = (\infty,\dots,\infty)$ for all $x\in\partial\Lambda$.
\end{lemma}

\begin{proof}
For any $v\in\Lambda^0$, compact neighbourhood $U\subset\Lambda^0$ of $v$ and $i\in\{1,\dots,k\}$, we have $U\Lambda^{e_i}\in v\CE(\Lambda)$. Therefore, given $x\in\partial\Lambda$, $m\le d(x)$ and $i\in\{1,\dots,k\}$, it follows that $x(m,m+e_i)\in x(m)\Lambda^{e_i}$, which can only occur if $d(x)=(\infty,\dots,\infty)$.
\end{proof}

\begin{lemma}\label{lem:H is amenable}
Let $(\Lambda,d)$ be a proper topological $k$-graph. The groupoid $H_\Lambda = \Zz(\Lambda*_s\Lambda,0)$, comprising all elements in $\Gg_\Lambda$ of the form $(x,0,y)$, is amenable.
\end{lemma}

\begin{proof}
For $m\in\NN^k$, define
\begin{align*}
H(m) &:= \{(\lambda x,0,\mu x)\in\Gg_\Lambda ~|~ d(\lambda)=d(\mu) \le m\} \\
&= \bigcup_{n\le m}\Zz(\Lambda^n*_s\Lambda^n,0).
\end{align*}
Then each $H(m)$ is a subgroupoid of $\Gg_\Lambda$, each has $\Gg_\Lambda^{(0)}$ as its unit space, and for $m\le n$, $H(m)$ is both an open and closed subgroupoid of $H(n)$. Hence we have a direct system of groupoids $\{H(m) ~|~ m\in\NN^k\}$ with direct-limit groupoid $H_\Lambda := \bigcup_{m\in\NN^k}H(m) = \Zz(\Lambda*_s\Lambda,0)$.

We claim that for each $m\in\NN^k$, $H(m)$ is a proper groupoid: To see this, let $W\subset \Gg_\Lambda^{(0)}\times\Gg_\Lambda^{(0)}$ be compact. There exist compact sets $U_i,V_j\subset\Lambda^0$, for $i=1,\dots,l_1$ and $j=1,\dots,l_2$, such that
\[
W\subset \bigcup_{i=1}^{l_1}\bigcup_{j=1}^{l_2}\Zz(U_i)\times \Zz(V_j).
\]
We then have
\begin{align}
(r,s)^{-1}&\Big( \bigcup_{i=1}^{l_1}\bigcup_{j=1}^{l_2}\Zz(U_i)\times \Zz(V_j)\Big) \label{eqn:inverse image}\\
&= \Big\{ (x,0,y)\in H(m) ~|~ x\in \bigcup_{i=1}^{l_1}\Zz(U_i), y\in \bigcup_{j=1}^{l_2}\Zz(V_j)\Big\} \notag\\
&= \bigcup_{n\le m}\bigcup_{i=1}^{l_1}\bigcup_{j=1}^{l_2}\Zz(U_i\Lambda^n*_sV_j\Lambda^n, 0), \notag
\end{align}
which is compact since the $U_i\Lambda^n, V_j\Lambda^n$ are compact on account of $\Lambda$ being proper. Since $(r,s)^{-1}(W)$ is a closed subset of \eqref{eqn:inverse image}, it follows that $(r,s)^{-1}(W)$ is compact. Thus each $H(m)$ is a proper groupoid, and Proposition~\ref{prop:proper implies amenable} implies each $H(m)$ is amenable. By \cite[Proposition~5.3.37]{A-DRen}, the direct limit $H_\Lambda$ is amenable.
\end{proof}

\begin{theorem}
Let $(\Lambda,d)$ be a proper topological higher-rank graph without sources. Then $\Gg_\Lambda$ is amenable.
\end{theorem}

\begin{proof}
Let $c:\Gg_\Lambda\to\ZZ^k$ be the continuous functor given by $c(x,m,y) = m$. We will show that the skew-product groupoid $\Gg_\Lambda(c)$ is amenable; the result will then follow from \cite[Proposition~II.3.8]{Ren1}.

We identify the unit space $\Gg_\Lambda(c)^{(0)}$ with $\Gg_\Lambda^{(0)}\times\ZZ^k$, and for each $m\in\ZZ^k$, define $U_m := \Gg_\Lambda^{(0)}\times\{m\}$. Each $\Gg_\Lambda(c)|_{U_m}$ is isomorphic to $H_\Lambda$, so by Lemma~\ref{lem:H is amenable}, each $\Gg_\Lambda(c)|_{U_m}$ is amenable.

For $m\in\ZZ^k$, define $H_m := \Gg_\Lambda(c)|_{[U_m]}$, where $[U_m] := s(r^{-1}(U_m))$ is the saturation of $U_m$. By \cite[Example~2.7]{MRW}, $\Gg_\Lambda(c)|_{U_m}$ is equivalent to $H_m$, so it follows from \cite[Theorem~2.2.17]{A-DRen} that each $H_m$ is amenable. We also have $[U_m]\subset [U_n]$ whenever $m\le n$, so, defining a cofinal sequence $\langle m_j\rangle_{j\in\NN}$ in $\NN^k$, we have $\bigcup_{j\in\NN} [U_{m_j}] = \Gg_\Lambda(c)^{(0)}$.

For each $j\in\NN$, amenability gives $C^*(H_{m_j})=C^*_{\rm red}(H_{m_j})$, and since $[U_{m_j}]$ is open in $\Gg_\Lambda(c)^{(0)}$, there exists a homomorphism $\pi_j:C^*(H_{m_j})\to C^*_{\rm red}(\Gg_\Lambda(c))$ defined by the inclusion $C_c(H_{m_j})\to C_c(\Gg_\Lambda(c))$. Furthermore, amenability implies that each $C^*(H_{m_j})$ is nuclear, so it follows that each image $\pi_j(C^*(H_{m_j}))$ is nuclear. We also have $\pi_j(C^*(H_{m_j}))\subset \pi_j(C^*(H_{m_{j+1}}))$ for each $j\in\NN$, and $\overline{\bigcup_{j\in\NN}\pi_j(H_{m_j})} = C^*_{\rm red}(\Gg_\Lambda(c))$. Therefore \cite[Theorem~2.3.9]{Lin} implies $C^*_{\rm red}(\Gg_\Lambda(c))$ is nuclear, and it follows from \cite[Corollary~6.2.14(ii) and Theorem~3.3.7]{A-DRen} that $\Gg_\Lambda(c)$ is amenable. Since $\ZZ^k$ is amenable, it follows from \cite[Proposition~II.3.8]{Ren1} that $\Gg_\Lambda$ is amenable.
\end{proof}

\section{$C^*$-algebras of topological higher-rank graphs}

\begin{examples}
\begin{enumerate}
\item
Let $(\Lambda,d)$ be a finitely aligned discrete $k$-graph. It follows from Examples~\ref{ex:path-groupoid examples}\eqref{ex:item:k-graph path groupoid} and~\ref{ex:boundary-path groupoids}\eqref{ex:item:k-graph boundary-path groupoid} together with \cite[Theorem~6.9 and~Theorem~6.13]{FMY1} that $C^*(G_\Lambda)\cong\Tt C^*(\Lambda)$ and $C^*(\Gg_\Lambda)\cong C^*(\Lambda)$, where $\Tt C^*(\Lambda)$ and $C^*(\Lambda)$ are defined and studied in \cite{RS} and \cite{RSY2}, respectively. Note that this example includes the Toeplitz and Cuntz-Krieger algebras of arbitrary directed graphs as studied in \cite{FowR2,BPRSz,FowLR,Pat2,RSz,DrTom2,R3} (among others).
\item
As we saw in Example~\ref{ex:graph examples}\eqref{ex:item:top graph}, there is a one-to-one correspondence between topological $1$-graphs and second-countable topological graphs. Given a topological $1$-graph $(\Lambda,d)$ with corresponding topological graph $E_\Lambda$, \cite[Theorems~5.1 and~5.2]{Y2} say that $C^*(G_\Lambda)\cong \Tt(E_\Lambda)$ and $C^*(\Gg_\Lambda)\cong \Oo(E_\Lambda)$, where $\Tt(E_\Lambda)$ and $\Oo(E_\Lambda)$ are, respectively, the Toeplitz and Cuntz-Krieger algebras of the topological graph $E_\Lambda$, as defined in \cite{Kat1}.
\item
Let $X$ be a second-countable locally compact Hausdorff space and let $\{\theta_i\}_{i=1}^k$ be a family of commuting homeomorphisms of $X$ onto itself. There is an induced action $\Theta$ of $\ZZ^k$ on $C_0(X)$ defined by
\[
\Theta_m(f)(x) = f(\theta_1^{m_1}\cdots\theta_k^{m_k}(x)),
\]
with universal crossed product $(C_0(X)\times_\Theta\ZZ^k,j_{C_0(X)},j_{\ZZ^k})$.

Recalling the topological $k$-graph $(\Lambda(X,\{\theta\}_{i=1}^k),d)$ defined in Example~\ref{ex:graph examples}\eqref{ex:item:k-GDS}, we have $C^*(\Gg_{\Lambda(X,\{\theta_i\})})\cong C_0(X)\times_\Theta\ZZ^k$; there are a number of ways to see this: for example, one can show there is a covariant representation $(\rho_{C_0(X)},\rho_{\ZZ^k})$ of the dynamical system $(C_0(X),\Theta,\ZZ^k)$ in $C^*(\Gg_{\Lambda(X,\{\theta_i\})})$ which induces an isomorphism $\rho_{C_0(X)}\times\rho_{\ZZ^k}: C_0(X)\times_\Theta\ZZ^k\to C^*(\Gg_{\Lambda(X,\{\theta_i\})})$; alternatively, one can show that $\Gg_{\Lambda(X,\{\theta_i\})}$ is isomorphic to the transformation groupoid $X\times\ZZ^k$, whose $C^*$-algebra is, in turn, isomorphic to the crossed-product $C^*$-algebra.
\end{enumerate}
\end{examples}

The above examples allow us to coherently make the following definition.

\begin{dfn}
Let $(\Lambda,d)$ be a compactly aligned topological $k$-graph. We define the \emph{Toeplitz algebra of $\Lambda$} to be the full groupoid $C^*$-algebra $C^*(G_\Lambda)$, and we define the \emph{Cuntz-Krieger algebra of $\Lambda$} to be the full groupoid $C^*$-algebra $C^*(\Gg_\Lambda)$.
\end{dfn}

\section{skew-product topological higher-rank graphs and crossed products by coactions}

In this section we extend the definition of a skew-product $k$-graph to topological $k$-graphs, and show that the associated groupoids can be realized as skew-product groupoids, extending \cite[Theorem~5.2]{KP2}. We can then realize certain crossed product $C^*$-algebras as topological higher-rank graph $C^*$-algebras, building on \cite[Corollary~5.3]{KP2} and \cite[Theorem~2.4]{KalQR}.

Given a topological groupoid $\Gamma$, a locally compact group $A$ and a continuous functor $b:\Gamma\to A$, we denote the \emph{skew product of $\Gamma$ by $b$} as $\Gamma(b)$; that is, $\Gamma(b)$ is the locally compact groupoid obtained by defining on $\Gamma\times A$ the multiplication $(x,g)(y,gb(x)):=(xy,g)$ and the inverse $(x,g)^{-1}:= (x^{-1},gb(x))$. In our setting, the groupoid $\Gamma$ is $r$-discrete and admits a Haar system, and it follows that the same is true for $\Gamma(b)$.

\begin{dfn}
Let $(\Lambda,d)$ be a topological $k$-graph, let $A$ be a locally compact group, and let $c:\Lambda\to A$ be a continuous functor. Define $\Lambda\times_c A$ to be the category with object and morphism sets
\[
\Obj(\Lambda\times_cA) = \Obj(\Lambda)\times A \quad\text{and}\quad \Mor(\Lambda\times_cA) = \Mor(\Lambda)\times A,
\]
range and source maps
\[
r(\lambda,a)=(r(\lambda),a) \quad\text{and}\quad s(\lambda,a) = (s(\lambda),ac(\lambda)),
\]
and composition
\[
(\lambda,a)(\mu,ac(\lambda)) = (\lambda\mu,a).
\]
Define a functor $d:\Lambda\times_cA\to\NN^k$ by $d(\lambda,a) = d(\lambda)$. Then, giving the object and morphism sets their product topologies, the pair $(\Lambda\times_cA,d)$ is a topological $k$-graph, called the \emph{skew-product of $(\Lambda,d)$ by $c$}.
\end{dfn}

\begin{lemma}
If $(\Lambda,d)$ is compactly aligned, then so is $(\Lambda\times_cA,d)$.
\end{lemma}

\begin{proof}
Let $p,q\in\NN^k$, and let $E\subset(\Lambda\times_cA)^p$ and $F\subset(\Lambda\times_cA)^q$ be compact. Let $P_\Lambda:\Lambda\times_cA\to\Lambda$ and $P_A:\Lambda\times_cA\to A$ be the coordinate maps $P_\Lambda(\lambda,a)=\lambda$ and $P_A(\lambda,a)=a$. We then see that $E\vee F$ is compact since it is a closed subset of the compact set $(P_\Lambda(E)\vee P_\Lambda(F))\times (P_A(E)\cap P_A(F))$.
\end{proof}

The proof of the following lemma is straightforward.

\begin{lemma}\label{lem:functor on Lambda induces one on G}
Let $(\Lambda,d)$ be a compactly aligned topological $k$-graph, let $A$ be a locally compact topological group, and let $c:\Lambda\to A$ be a continuous functor. Then there is a continuous functor $\tilde{c}:G_\Lambda\to A$ defined by
\[
\tilde{c}(\lambda x,d(\lambda)-d(\mu),\mu x) = c(\lambda)c(\mu)^{-1}.
\]
\end{lemma}

\begin{prop}\label{prop:skew products match}
Let $(\Lambda,d)$ be a compactly aligned topological $k$-graph, let $A$ be a locally compact topological group, and let $c:\Lambda\to A$ be a continuous functor. Then, with the notation of Lemma~\ref{lem:functor on Lambda induces one on G},
\[
G_\Lambda(\tilde{c})\cong G_{\Lambda\times_cA}.
\]
Furthermore, denoting the restriction of $\tilde{c}$ to $\Gg_\Lambda = G_\Lambda|_{\partial\Lambda}$ again by $\tilde{c}$, we have
\[
\Gg_\Lambda(\tilde{c})\cong \Gg_{\Lambda\times_cA}.
\]
\end{prop}

\begin{proof}
We first define a functor $\phi:G_\Lambda(\tilde{c})\to G_{\Lambda\times_cA}$. For $(x,a)\in (G_\Lambda(\tilde{c}))^{(0)}$, define a path $\phi(x,a):\Omega_{k,d(x)}\to\Lambda\times_cA$ in $X_{\Lambda\times_cA}$ by
\[
(\phi(x,a))(m,n) = (x(m,n),ac(x(0,m))) \quad\text{for } m\le n\le d(x).
\]
Now let $((\lambda x,d(\lambda)-d(\mu),\mu x),a)\in G_\Lambda(\tilde{c})$ and define
\begin{align*}
\phi((\lambda x,d(\lambda)-d(\mu),\mu x),a) &= (\phi(\lambda x,a),d(\lambda)-d(\mu),\phi(\mu x, a\tilde{c}(\lambda x,d(\lambda)-d(\mu),\mu x))) \\
&= (\phi(\lambda x,a),d(\lambda)-d(\mu),\phi(\mu x, ac(\lambda)c(\mu)^{-1})).
\end{align*}
Straightforward but lengthy calculations then show that $\phi:G_\Lambda(\tilde{c})\to G_{\Lambda\times_cA}$ is a bijective continuous functor with continuous inverse, and the first part of the proposition follows.

We now show that $\Gg_\Lambda(\tilde{c})\cong \Gg_{\Lambda\times_c A}$. By definition, we have $\Gg_\Lambda(\tilde{c})=G_\Lambda(\tilde{c})|_{(\partial\Lambda)\times A}$, so it suffices to show that $\phi((\partial\Lambda)\times A) = \partial(\Lambda\times_c A)$.

First fix $\phi(x,a)\in \partial(\Lambda\times_c A)$; we show that $(x,a)\in (\partial\Lambda)\times A$. Let $m\in\NN^k$ satisfy $m\le d(x)$, and let $E\in x(m)\CE(\Lambda)$. Choosing any compact neighbourhood $B$ of $ac(x(0,m))$, we have $E\times B\in (\phi(x,a)(m))\CE(\Lambda\times_c A)$, so there exists $(\lambda,b)\in E\times B$ such that $\phi(x,a)(m,m+d(\lambda)) = (\lambda,b)$. We then have $x(m,m+d(\lambda))=\lambda\in E$, giving $\partial(\Lambda\times_c A)\subset \phi((\partial\Lambda)\times A)$.

On the other hand, fix $(x,a)\in (\partial\Lambda)\times A$, let $m\in \NN^k$ satisfy $m\le d(\phi(x,a))=d(x)$, and let $E\in (\phi(x,a)(m))\CE(\Lambda\times_c A)$.

Let $\{\Bb_j\}_{j\in\NN}$ be a neighbourhood basis for $\phi(x,a)(m) = (x(m),ac(x(0,m)))$ such that each $\Bb_j\subset r(E)$ and $\overline{\Bb_{j+1}}\subset \Bb_j$. For each $j\in\NN$, define
\[
F_j := \overline{\Bb_j}(\Lambda\times_c A)\cap E = \{(\lambda,b)\in E : r(\lambda,b)\in\overline{\Bb_j}\},
\]
so $F_j\in (\phi(x,a)(m))\CE(\Lambda\times_c A)$.

Let $P_\Lambda:\Lambda\times_c A\to\Lambda$ denote the coordinate map. For each $j\in\NN$, we have $P_\Lambda(F_j)\in x(m)\CE(\Lambda)$, so there exists $(\lambda_j,b_j)\in F_j$ such that $x(m,m+d(\lambda_j))=\lambda_j$. Since each $r(\lambda_j,b_j)\in\Bb_j$, it follows that
\begin{equation}\label{eqn:range of sequence}
\lim_{j\in\NN} (r(\lambda_j),b_j) = (x(m),ac(x(0,m))).
\end{equation}

Since $\langle (\lambda_j,b_j)\rangle_{j\in\NN}$ is contained in the compact set $E$, there exists a convergent subsequence $\langle (\lambda_j,b_j)\rangle_{j\in J}$ with limit $(\lambda,b)\in E$. We then have
\begin{align*}
(r(\lambda),b) = r(\lambda,b) &= r\big(\lim_{j\in J}(\lambda_j,b_j)\big) \\
&= \lim_{j\in\NN}r(\lambda_j,b_j) = (x(m),ac(x(0,m))) \quad\text{by \eqref{eqn:range of sequence}},
\end{align*}
hence $b=ac(x(0,m))$. Furthermore,
\[
\lambda = \lim_{j\in J}\lambda_j=\lim_{j\in J}x(m,m+d(\lambda_j))=x(m,m+d(\lambda)),
\]
and it follows that $\phi(x,a)(m,m+d(\lambda))=(\lambda,b)\in E$, and $\phi(x,a)\in\partial(\Lambda\times_c A)$.

Therefore $\phi((\partial\Lambda)\times A) = \partial(\Lambda\times_c A)$, so $\phi$ restricts to an isomorphism from $\Gg_\Lambda(\tilde{c})$ onto $\Gg_{\Lambda\times_c A}$, completing the proof.
\end{proof}

\begin{notation}\label{notation:actions}
For a continuous functor $b$ from a locally compact groupoid $G$ with continuous Haar system to a locally compact abelian group $A$, we denote by $\alpha(b)$ the action of the dual group $\widehat A$ on $C^*(G)$ defined by \cite[Proposition~II.5.1]{Ren1}.

On the other hand, for a continuous functor $b$ from a locally compact $r$-discrete groupoid $G$ with Haar system to a discrete group $A$, we denote by $\delta(b)$ the coaction of $A$ on $C^*(G)$ defined by \cite[Lemma~4.2]{KalQR}.
\end{notation}

\begin{rmk}
Our Theorem~\ref{thm:coaction of group} concerns coactions of discrete groups (see \cite{Ng,Q2,KalQR,R3}, for example). There is, however, much literature on the theory of coactions of locally compact groups on $C^*$-algebras (see \cite{LanPhilRSuth,R1,Q1,R2}, for example). Our reliance on groupoid theory in the proofs of Theorem~\ref{thm:action of dual group} and Theorem~\ref{thm:coaction of group} has meant that our theorems only address actions of locally compact abelian groups and coactions of discrete groups. It is possible that Theorem~\ref{thm:coaction of group} holds in the generality of locally compact groups; consequently, Theorem~\ref{thm:action of dual group} would follow as the abelian case.
\end{rmk}

\begin{theorem}\label{thm:action of dual group}
Let $(\Lambda,d)$ be a compactly aligned topological $k$-graph, let $A$ be a locally compact abelian group, and let $c:\Lambda\to A$ be a continuous functor. Then, with the notation of Lemma~\ref{lem:functor on Lambda induces one on G} and Notation~\ref{notation:actions},
\[
C^*(G_\Lambda)\times_{\alpha(\tilde{c})}\widehat{A} \cong C^*(G_{\Lambda\times_cA}).
\]
Furthermore, denoting the restriction of $\tilde{c}$ to $\Gg_\Lambda$ again by $\tilde{c}$, we have
\[
C^*(\Gg_\Lambda)\times_{\alpha(\tilde{c})}\widehat{A} \cong C^*(\Gg_{\Lambda\times_cA}).
\]
\end{theorem}

\begin{proof}
Both parts of of the theorem are achieved in two steps, using \cite[Theorem~II.5.7]{Ren1} and Proposition~\ref{prop:skew products match}.
\end{proof}

\begin{theorem}\label{thm:coaction of group}
Let $(\Lambda,d)$ be a compactly aligned topological $k$-graph, let $A$ be a discrete group and let $c:\Lambda\to A$ be a continuous functor. Then, with the notation of Lemma~\ref{lem:functor on Lambda induces one on G} and Notation~\ref{notation:actions},
\[
C^*(G_\Lambda)\times_{\delta(\tilde{c})} A\cong C^*(G_{\Lambda\times_cA}).
\]
Furthermore, denoting the restriction of $\tilde{c}$ to $\Gg_\Lambda$ again by $\tilde{c}$, we have
\[
C^*(\Gg_\Lambda)\times_{\delta(\tilde{c})}A\cong C^*(\Gg_{\Lambda\times_cA}).
\]
\end{theorem}

\begin{proof}
The theorem follows from \cite[Theorem~4.3 and Theorem~6.2]{KalQR} and Proposition~\ref{prop:skew products match}.
\end{proof}

\end{document}